\documentclass[journal]{IEEEtran}

\usepackage{amsmath,amssymb,amsfonts,mathtools}
\usepackage{graphicx}
\usepackage{booktabs}
\usepackage[colorlinks,citecolor=gray,linkcolor=blue]{hyperref}
\usepackage{amsthm}
\usepackage{bbm}
\usepackage{bm} 
\usepackage{xcolor}
\usepackage{tikz}
\usetikzlibrary{automata,positioning,arrows.meta}
\usepackage{pgfplots}
\pgfplotsset{compat=1.18}
\usepackage{bbold}   
\usepackage[numbers,sort&compress]{natbib}

\usepackage[T1]{fontenc}
\usepackage[utf8]{inputenc}

\usepackage{amsthm}

\newtheorem{lemma}{Lemma}

\usepackage{times}
\usepackage{thmtools}
\usepackage[capitalise,noabbrev]{cleveref}
\usepackage{booktabs} 
\usepackage{multirow} 

\usepackage{bbm}
\usepackage{bm} 
\usepackage{graphicx}
\usepackage{amsmath}
\usepackage{amsthm}

\DeclareMathOperator*{\maximize}{maximize}
\DeclareMathOperator*{\minimize}{minimize}
\DeclareMathOperator*{\stc}{s.\,t.}

\theoremstyle{plain}
\newtheorem{theorem}{Theorem}
\newtheorem{proposition}{Proposition}
\theoremstyle{definition}
\newtheorem{definition}{Definition}
\newtheorem{example}{Example}
\theoremstyle{remark}

\newtheorem{assumption}{Assumption}

\newcommand{\tr}{^{\top}}

\newcommand{\Real}{\mathbb{R}}
\newcommand{\Nats}{\mathbb{N}_0}

\definecolor{NavyBlue}{rgb}{0.0, 0.0, 0.5} 

\renewcommand{\P}{\mathbb{P}}
\newcommand{\E}{\mathbb{E}}

\newcommand{\PiSR}{\Pi_{\mathrm{SR}}}
\newcommand{\PiHR}{\Pi_{\mathrm{HR}}}
\newcommand{\PiMR}{\Pi_{\mathrm{MR}}}

\newcommand{\opt}{^{\star}}

\begin{document}

\title{Computing Optimal Non-stationary Policies for the Probabilistic Reach--Avoid Problems}

\author{Saber Omidi, Momotaz Begum, Marek Petrik, Se Young Yoon
}
\maketitle

\begin{abstract}
We develop and analyze an algorithm for computing optimal non-stationary policies for probabilistic reach--avoid problems. The objective in probabilistic reach--avoid problems is to maximize the probability of reaching a goal without visiting predefined avoid states. Prior work has focused on computing optimal \emph{stationary} reach--avoid policies. Stationary policies use the same action distribution
whenever a given state is visited.This paper shows that stationary reach--avoid policies may
be suboptimal and non-stationary policies may outperform
them in reach probability by a margin arbitrarily close to one. To compute optimal non-stationary policies, we reduce the reach--avoid problem to a constrained average-reward Markov decision problem and propose an algorithm to compute lottery-based policies from a standard linear program formulation. We prove that the lottery-based policies are optimal and show that they are simple to compute, analyze, interpret, and implement. The numerical experiments confirm that the lottery-based policies can outperform state-of-the-art stationary policies even in complex domains. 
\end{abstract}

\begin{IEEEkeywords}
Probabilistic reach--avoid control, constrained average-reward Markov decision processes, linear programming, randomized switching policies, probabilistic safe sets, safety-critical control.
\end{IEEEkeywords}

\section{Introduction}\label{sec:intro}

Reach--avoid problems seek a controller that steers a system
to a designated set of \emph{reach state} while avoiding \emph{avoid states},
such as those representing obstacles or hazardous regions. Probabilistic reach--avoid systems also allow environments with significant uncertainty that produce stochastic dynamics~\cite{xue2023reach}. Such reach--avoid models apply to a wide range of applications, including autonomous vessels, mobile robots, and many others~\cite{pietrzykowski2022autonomous,fan2021controller}. In these applications, the controller must reach a predetermined goal while avoiding states that could lead to damage, injury, or other significant loss. 

In this paper, we study the infinite-horizon, probabilistic reach--avoid problem, where we maximize the probability of visiting a reach state while keeping the probability of visiting an avoid state sufficiently small. Following prior literature, we model the problem using the standard discrete-time control framework of Markov decision processes~(MDPs)~\cite{Puterman2005,kallenberg2011markov}. In particular, {\'A}vila and Junca~\cite{avila2021reachability} show that probabilistic reach--avoid problems can be formulated as MDPs with the standard average-reward objective. 
The average-reward formulation requires careful analysis
of long-run state--action frequencies.
In constrained multichain MDPs, a stationary policy need not
realize the frequencies prescribed by an optimal linear-program
solution, so policy construction requires separate
justification~\cite{kallenberg1983linear}.

Before stating our contributions, we briefly summarize the average-reward formulation of {\'A}vila and Junca~\cite{avila2021reachability}, which we build on. As stated above, their pioneering work shows that the reach--avoid objective can be analogously expressed as the long-run average of an MDP with an appropriately constructed reward. The formulation (P2) in~\cite{avila2021reachability} is equivalent to our reach--avoid objective with only superficial differences. In their formulation, they assign a reward of 1 when the controller visits a reach state and 0 otherwise. To keep the probability of visiting an avoid state sufficiently small, they use a Lagrangian relaxation to assign a flexible penalty for visiting an avoid set, then optimize over the penalty. For a fixed penalty value, they use a linear program to compute the optimal average-reward value function, commonly referred to as \emph{gain}. Finally, they use the optimal gain to construct a stationary policy using the standard construction in average-reward MDPs~\cite[theorem~9.3.8]{Puterman2005}. 

It is important to note that the policies of {\'A}vila and Junca~\cite{avila2021reachability} are stationary only after the states of the model are augmented with an additional binary variable that indicates whether the process has previously visited one of the avoid states. We emphasize that the need for this additional indicator variable depends on whether the controller is allowed to continue after visiting an avoid state or must stop. In the remainder of this discussion, we assume that the state space is appropriately augmented and the policies of {\'A}vila and Junca~\cite{avila2021reachability} are stationary in this space. 
Stationarity therefore depends on the chosen state space.

We make the following main contributions. First, we formulate the probabilistic reach--avoid problem
as a constrained average-reward objective in MDPs using
the limiting expected state--action frequencies. Our constrained formulation builds on limiting distributions instead of gain and simplifies the derivation from prior work~\cite{avila2021reachability}. Specifically, our derivation does not require Lagrangian relaxation to formulate a linear program.

Second, building on our constrained linear program formulation, we prove that stationary policies are insufficient to optimize the probabilistic reach--avoid objective. We give an example in which the  suboptimality gap of the reach probability between any stationary and best Markov policy is arbitrarily close to 1.  This result indicates that the policies proposed by {\'A}vila and Junca~\cite{avila2021reachability} are inherently suboptimal. Our counterexample reflects the critical difference between constrained and unconstrained multichain average-reward MDPs~\cite{hordijk1984constrained,krass2002achieving,beutler1985optimal,ross1989randomized,altman1991markov}. Our work further investigates the classical work on randomized and past-dependent policies~\cite{beutler1985optimal,ross1989randomized} and state--action frequencies~\cite{altman1991markov} in constrained average-reward MDPs.

Third, we adapt a construction of Krass and Vrieze~\cite{krass2002achieving} to propose an algorithm for constructing optimal non-stationary (history-dependent or equivalent Markov) policies for the reach--avoid problem. Our algorithm begins by solving the standard constrained linear-program formulation of average-reward MDPs. This formulation incorporates two separate reward functions: one for reach states and another one for the constrained states. Then, the algorithm constructs a lottery-based policy from the linear program's optimal solution. We note that even though our derivation of the linear program follows very different steps, we recover a linear program that is identical to {\'A}vila and Junca~\cite{avila2021reachability}. The main distinction, however, is that the proposed lottery-based policies are guaranteed to achieve the optimal value of the linear program whenever it is feasible. In contrast, for some problem instances, no stationary policy that satisfies the safety constraint can achieve this optimal value.

The rest of the paper is structured as follows. 
\Cref{sec:ra_formulation} introduces the probabilistic
reach--avoid problem and shows that stationary policies
may fail to maximize the reach probability subject to
the prescribed safety constraint.
\Cref{sec:avg} provides the constrained average-reward MDP formulation, specifying the MDP structure and the average-reward criteria for reachability and safety. \Cref{sec:CLP} presents the linear programming formulation and constructs optimal non-stationary lottery-based policies from optimal LP solutions, using a single random switch between two stationary policies. \Cref{sec:rp} characterizes probabilistic safe sets using the optimal safety gain, compares the proposed reach--avoid formulation with the formulation of {\'A}vila and Junca, and discusses discounted approximations of the almost-sure safe set. Finally, \cref{sec:nv} evaluates the proposed algorithm numerically by comparing probabilistic safe sets with discounted approximations for the inverted pendulum and comparing lottery-based policies with stationary policies for Zermelo's navigation problem.

\section{Probabilistic Reach--Avoid Formulation}
\label{sec:ra_formulation}

In this section, we introduce our notation and define our objective formally.  The symbols $\mathbb{R}$, $\mathbb{R}_{\geq 0}$, and $\mathbb{N}_0$ refer to reals, nonnegative reals, and nonnegative integers, respectively. For a finite set $\mathcal{X}$, $\Delta(\mathcal{X})$ denotes the probability simplex and $\mathbb{1}_{\mathcal{X}} \colon \mathcal{Y} \to  \{0,1\} $ denotes the indicator function of $\mathcal{X} \subseteq \mathcal{Y}$.  Vectors are denoted by lowercase bold symbols, such as $\bm{v} \in \Real^n$, with $v_i \in \Real$ representing the $i$-th entry. Matrices are denoted by uppercase bold symbols, such as $\bm{P} \in \Real^{n \times m}$, where $P_{ij}$ denotes the $(i,j)$-th entry, $\bm{I}$ denotes the identity matrix of dimension appropriate to the context.

As in prior work, we model the underlying system as a finite \emph{controlled Markov chain} $\left(\mathcal{S}, \{\mathcal{A}(s)\}_{s\in\mathcal{S}}, p, \alpha \right)$~\cite{avila2021reachability}. Here, $\mathcal{S}$ is the finite non-empty \emph{state set}, and $\mathcal{A}(s),\, s \in \mathcal{S}$ are the finite non-empty \emph{action sets}.  To streamline our notation, we use $\mathcal{Z} := \{(s,a) \mid  s\in\mathcal{S}, \,a\in\mathcal{A}(s)\}$ to denote admissible \emph{state-action} pairs and $\mathcal{A} := \bigcup_{s\in \mathcal{S}} \mathcal{A}(s)$ to denote the set of \emph{all actions}. Then, $p\colon \mathcal{Z} \to \Delta(\mathcal{S})$ is the transition probability function and $\alpha\in\Delta(\mathcal{S})$ is the \emph{initial distribution}. The notation $p(s' \mid s, a)$ represents the probability of transitioning from $s\in \mathcal{S}$ to $s' \in \mathcal{S}$ after taking $a\in \mathcal{A}(s)$.

A crucial concept in formalizing the reach--avoid problem is to define a \emph{history}, also referred to as a \emph{trajectory} in reinforcement learning, that represents a sample of the execution of the process. The set of \emph{histories} of length $t \in \Nats$ is
\[
\mathcal{H}_t := \left\{(s_0,a_0,\ldots,s_t) \mid s_i \in \mathcal{S}, a_i \in \mathcal{A}(s_i), \forall i \right\}.
\]
Let $\mathcal{H}_{\infty} = \cup_{t\in \Nats} \mathcal{H}_t$. Below, we use $\mathcal{H}_{\infty}$ as a sample space of a probability space and define, for each history $h = (s_0, a_0, \dots, s_t) \in \mathcal{H}_t$, a state random variable $\tilde{s}_k(h) := s_k, k = 0, \dots, t$ and an action random variable $\tilde{a}_k(h) := a_k, k = 0, \dots , t-1$. For any set $\mathcal{X} \subseteq \mathcal{S}$, we define the \emph{first hitting time} $\tilde{\tau}_{\mathcal{X}} \colon \mathcal{H}_{\infty} \to \Nats\cup\{\infty\}$ as
\begin{equation}
 \label{eq:hitting_times}
    \tilde{\tau}_{\mathcal{X}} := \min\left\{k\in\Nats \mid  \tilde{s}_k\in\mathcal{X} \right\},
\end{equation}
with the convention $\inf\emptyset=\infty$.

A policy for a controlled Markov chain is a sequence of decision rules that map finite histories of states and actions to probability distributions over admissible actions~\cite{puterman1990markov}.
The set of \emph{history-dependent randomized policies} is 
\begin{equation*}
\begin{aligned}
  \PiHR := \Bigl\{
  &(\pi_0,\pi_1,\ldots) \mid \pi_i \colon \mathcal H_i \to \Delta(\mathcal A),\ i\in\Nats,\\
  &\pi_i(a\mid h)=0,\ \forall h=(\ldots,s_i)\in\mathcal H_i,\ a\notin\mathcal A(s_i)\Bigr\}.
\end{aligned}
\end{equation*}
Here, $\pi_i(a \mid h)$ is the probability of taking action $a$ given the finite history $h\in\mathcal H_i$, and $s_i$ is the final state of the history $h$. In \emph{Markov randomized policies}, defined as 
\begin{equation*}
\begin{aligned}
  \PiMR := \Bigl\{
  &(\pi_0, \pi_1, \dots ) \in \PiHR \mid \pi_i(a \mid h_1) = \pi_i(a \mid h_2), \\
&\forall h_1,h_2:\ \tilde{s}_i(h_1)=\tilde{s}_i(h_2),\ i\in\Nats \Bigr\},
\end{aligned}
\end{equation*}
each decision depends only on the last state of the history. It is well known that Markov policies can be identified with a sequence of simple decision rules: 
\begin{equation*}
\begin{aligned}
\PiMR = \Bigl\{ &\pi = (\pi_0, \pi_1, \dots ) \mid  \pi_i \colon \mathcal{S} \to \Delta(\mathcal{A}), \\
&\pi_i(a \mid s) = 0,\ \forall s\in \mathcal{S},\
a \notin \mathcal{A}(s),\ i\in \Nats \Bigr\}.
\end{aligned}
\end{equation*}
Finally, a policy is \emph{stationary randomized}  when the decision rule is time-invariant:
\[\PiSR := \left\{ (\pi_0, \pi_1, \dots  ) \in \PiMR \mid  \pi_i = \pi_j, \forall i,j \in \Nats \right\}.\] 
In most decision-theoretic applications, it is common to focus on stationary deterministic policies. A stationary policy $\pi \in \PiSR$ is \emph{deterministic} if its decision rule for each state is a degenerate probability distribution. In this work, we allow randomized policies because optimal
policies for constrained MDPs may need to be
randomized~\cite{altman2021constrained,beutler1985optimal,ross1989randomized}. Note that the need for non-stationary policies in probabilistic reach--avoid settings has not received much attention in the literature. 

Given a policy $\pi \in \PiHR$ and an initial distribution $\alpha$, we define the standard probability space over $\mathcal{H}_\infty=(\mathcal S\times\mathcal A)^{\Nats}$, equipped with the product $\sigma$-algebra, using the Ionescu--Tulcea theorem~(see, for example, \citep[appendix~B.2]{bauerle2011markov}). The probability measure and expectation induced by $\pi$ are denoted by $\mathbb{P}^{\pi}[\cdot]$ and $\mathbb{E}^{\pi}[\cdot]$. When the history starts from a state $s$, they are denoted by $\mathbb{P}^{\pi}_{s}[\cdot]$ and $\mathbb{E}^{\pi}_{s}[\cdot]$; when it starts from an initial distribution $\alpha$, they are denoted by $\mathbb{P}^{\pi}_{\alpha}[\cdot]$ and $\mathbb{E}^{\pi}_{\alpha}[\cdot]$.
\begin{definition}\label{def:RA}
The objective of the \emph{probabilistic reach--avoid} problem for a \emph{reach set} $\mathcal{R}\subseteq\mathcal{S}$ and \emph{avoid set} $\mathcal{O} \subseteq \mathcal{S} \setminus \mathcal{R}$ is defined as 
\begin{equation} \label{pr:RA}
\begin{aligned}
    \maximize_{\pi\in\PiHR}\quad & \mathbb{P}^{\pi}_{\alpha} \left[\tilde{\tau}_{\mathcal{R}} < \tilde{\tau}_{\mathcal{O}} \right] \\
    \stc\quad & \mathbb{P}^{\pi}_{\alpha} \left[\tilde{\tau}_{\mathcal{O}} < \infty \right] \le 1-\delta ,
\end{aligned}
\end{equation}
where $\delta\in [0,1]$ is the desired safety level. We also use the \emph{constraint set} $\mathcal{C} := \mathcal{S} \setminus \mathcal{O}$.
\end{definition}
Intuitively, the problem in~\eqref{pr:RA} seeks a policy that maximizes the probability of reaching a state in $\mathcal{R}$ before visiting a state in $\mathcal{O}$, while ensuring that the history never visits a state in $\mathcal{O}$ with probability at least $\delta$.

It will be useful to define an absorbing set as a set of states the process cannot leave after entering.
\begin{definition}
\label{def:absorbing_set}
A set $\mathcal B\subseteq\mathcal S$ is absorbing if
\begin{equation}
\label{eq:absorption}
\sum_{s'\in\mathcal B}p(s'\mid s,a)=1,
\qquad \forall s\in\mathcal B, \forall a\in \mathcal{A}(s).
\end{equation}
\end{definition}
We assume the following property in the rest of the paper.
\begin{assumption}
\label{ass:absorption}
The sets $\mathcal{R}$ and $\mathcal{O}$ are absorbing.
\end{assumption}

The assumption classifies trajectories according to whether they enter the reach or avoid set. Once a history enters either $\mathcal R$ or $\mathcal O$, it remains in that set permanently. This property guarantees that outcomes persist over time, allowing the probabilities in \cref{def:RA} to be represented  as long-run expected averages of the corresponding indicator functions. 
The gain–probability identities are established in \Cref{lemma:gain_prob},
and the cemetery-set interpretation of the absorbing avoid set
is discussed in \Cref{sec:rp}.

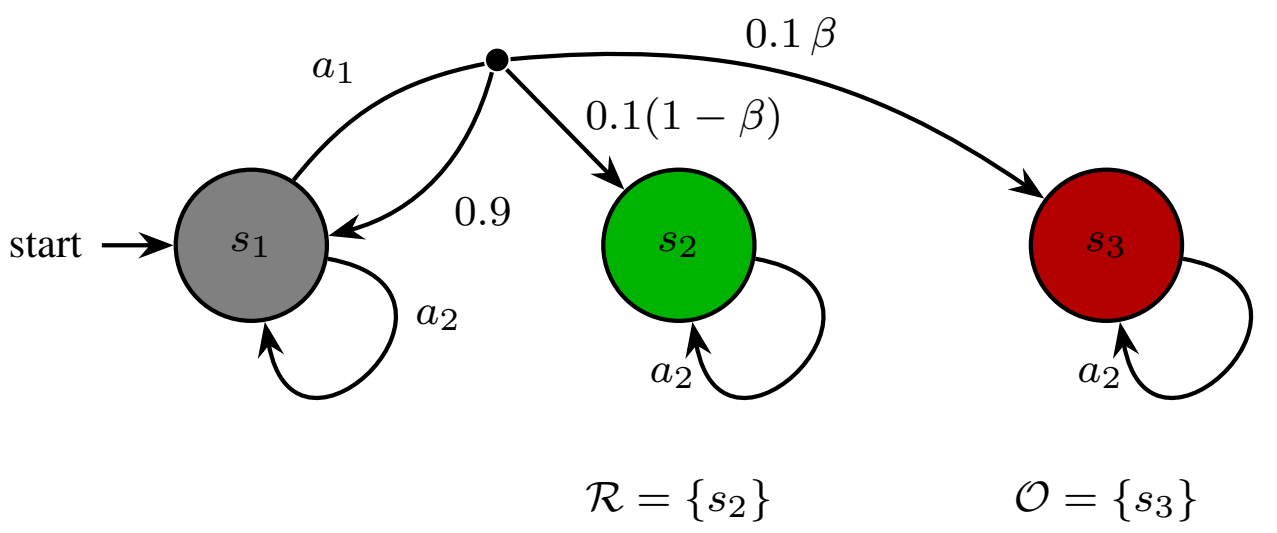
\begin{figure}
    \centering
    \begin{tikzpicture}[
        >=Stealth, thick, node distance=1.8cm, font=\footnotesize,
        every state/.style={draw, minimum size=1cm, inner sep=0pt},
        selfloop/.style={loop, out=-10, in=-80, looseness=5},
        lbl/.style={align=left},
        action/.style={circle, fill, inner sep=1.5pt},
        info/.style={below=9mm, align=center}]
      \node[state, initial, initial text=start, fill=gray] (s1) {$s_1$};
      \node[state, right=of s1, fill=green!70!black] (s2) {$s_2$};
      \node[state, right=of s2, fill=red!70!black] (s3) {$s_3$};
      \node[action, above right=8mm and 12mm of s1] (a1) {};
      \path[->]
        (s1) edge[-, bend left=20] node[lbl, above left] {$a_1$} (a1)
        (a1) edge[bend left=30] node[lbl, below right] {$0.9$} (s1)
             edge node[lbl, above right, yshift=-2mm] {$0.1(1-\beta)$} (s2)
             edge[bend left=20] node[lbl, above] {$0.1\,\beta$} (s3)
        (s1) edge[selfloop] node[lbl, right, pos=0.3] {$a_2$} (s1)
        (s2) edge[selfloop] node[lbl, left, pos=0.85] {$a_2$} (s2)
        (s3) edge[selfloop] node[lbl, left, pos=0.85] {$a_2$} (s3);
      \node[info] at (s2.south) {$\mathcal{R} = \left\{ s_2 \right\}$};
      \node[info] at (s3.south) {$\mathcal{O} = \left\{ s_3 \right\}$};
    \end{tikzpicture}
    \caption{Three-state reach--avoid example, parameterized by $\beta\in(0,1)$, illustrating the limitation of stationary policies in probabilistic reach--avoid problems.}
    \label{fig:example_toy}
\end{figure}

The following theorem shows that Markov policies can outperform stationary policies in reach probability by a margin arbitrarily close to one while achieving a desired safety level.
\begin{theorem}
\label{thm:stationary_insufficient}
For each $\varepsilon\in(0,1)$, there exists a Markov process, sets $\mathcal{R}, \mathcal{O} \subseteq \mathcal{S}$, a constant $\delta\in (0,1)$, and $\hat{\pi}\in\PiMR$ feasible in \Cref{pr:RA}, such that $\widehat{\Pi}_{\mathrm{SR}}\neq \emptyset$, and
\begin{equation*}
\mathbb P_\alpha^{\hat{\pi}}
[\tilde\tau_{\mathcal R}<\tilde\tau_{\mathcal O}]
\ge 
\mathbb P_\alpha^{\pi_{\mathrm{s}}}
[\tilde\tau_{\mathcal R}<\tilde\tau_{\mathcal O}] + \varepsilon,
\quad \forall \pi_{\mathrm{s}}\in\widehat{\Pi}_{\mathrm{SR}},
\end{equation*}
where $\widehat{\Pi}_{\mathrm{SR}} \subseteq \PiSR$ is the set of stationary policies feasible in the instance of~\Cref{pr:RA}.
\end{theorem}

\begin{example}\label{ex:toy}
Consider the controlled Markov chain in \cref{fig:example_toy} with
$\mathcal S=\{s_1,s_2,s_3\}$, $\alpha(s_1)=1$, $\mathcal R=\{s_2\}$, and
$\mathcal O=\{s_3\}$. The action sets are
$\mathcal A(s_1)=\{a_1,a_2\}$ and
$\mathcal A(s_2)=\mathcal A(s_3)=\{a_2\}$.
In $s_1$, action $a_2$ keeps the state at $s_1$, while
action $a_1$ transitions to $s_1$, $s_2$, and $s_3$ with probabilities
$p(s_1\mid s_1,a_1)=0.9$, $p(s_2\mid s_1,a_1)=0.1(1-\beta)$, and
$p(s_3\mid s_1,a_1)=0.1\,\beta$ for a parameter $\beta\in(0,1)$.
\Cref{ass:absorption} holds because states $s_2$ and $s_3$ are absorbing.
\end{example}

\begin{figure}
    \centering
    \begin{tikzpicture}
      \begin{axis}[
          width=0.8\columnwidth, height=5cm,
          xmin=-0.03, xmax=1.03, ymin=0.7, ymax=1.03,
          xtick={0,0.8,1},
          xticklabels={$0$, $1-\beta$, $1$},
          ytick={0.8,1},
          yticklabels={\rotatebox[origin=c]{90}{$1-\beta$}, $1$},
          xlabel={Reach: $\mathbb{P}^{\pi}_{s_1}[\tilde\tau_{\mathcal R}<\tilde\tau_{\mathcal O}]$},
          ylabel={Safety: $\mathbb{P}^{\pi}_{s_1}[\tilde\tau_{\mathcal O}=\infty]$},
          tick label style={font=\footnotesize},
          yticklabel style={inner sep=1pt},
          label style={font=\footnotesize},
          ylabel style={
            at={(axis description cs:-0.10,0.5)},
            anchor=south
          }
        ]
        \addplot[blue, very thick, dashed, domain=0:0.8]
          {1 - 0.25*x};

        \addplot[only marks, mark=*, mark size=3pt, red]
          coordinates {(0.8,0.8)};

        \addplot[only marks, mark=*, mark size=3pt, green!50!black]
          coordinates {(0,1)};

        \addplot[black, thick, dotted, domain=0:1] {0.84};

        \node[anchor=south east, font=\footnotesize, inner sep=1pt]
          at (axis cs:0.99,0.84) {$\delta$};

        \addplot[only marks, mark=*, mark size=2pt, black]
          coordinates {(0.64,0.84)};
        \node[
          anchor=north east,
          draw=black,
          fill=white,
          line width=0.3pt,
          inner sep=2pt,
          font=\fontsize{6}{7}\selectfont,
          align=left
        ] at (rel axis cs:0.99,0.99) {%
          \tikz[baseline=-0.5ex]{
            \draw[blue, very thick, dashed] (0pt,0pt)--(8pt,0pt);
          }\hspace{2pt}$\PiMR$\hspace{10pt}%
          \tikz[baseline=-0.5ex]{
            \fill[red] (0pt,0pt) circle[radius=1.6pt];
          }\hspace{2pt}$\pi_1\in\PiSR$\hspace{10pt}%
          \tikz[baseline=-0.5ex]{
            \fill[green!50!black] (0pt,0pt) circle[radius=1.6pt];
          }\hspace{2pt}$\pi_{\mathrm{s}}, \pi_2\in\PiSR$\\[2pt]
          \tikz[baseline=-0.5ex]{
            \draw[black, thick, dotted] (0pt,0pt)--(8pt,0pt);
          }\hspace{2pt}$\delta=1-\beta(1-\beta)$\hspace{10pt}%
          \tikz[baseline=-0.5ex]{
            \fill[black] (0pt,0pt) circle[radius=1.2pt];
          }\hspace{2pt}$\hat{\pi}_{\eta}\in\PiMR$
        };
      \end{axis}
    \end{tikzpicture}
    \caption{Reach probability and safety of stationary and Markov policies in \cref{ex:toy} as used in the proof of \cref{thm:stationary_insufficient}. Here, $\pi_1(a_1\mid s_1) = 1$ and $\pi_2(a_1\mid s_1) = 0$, and the dashed line represents the values achievable by $\PiMR$.}
    \label{fig:seg_line}
\end{figure}
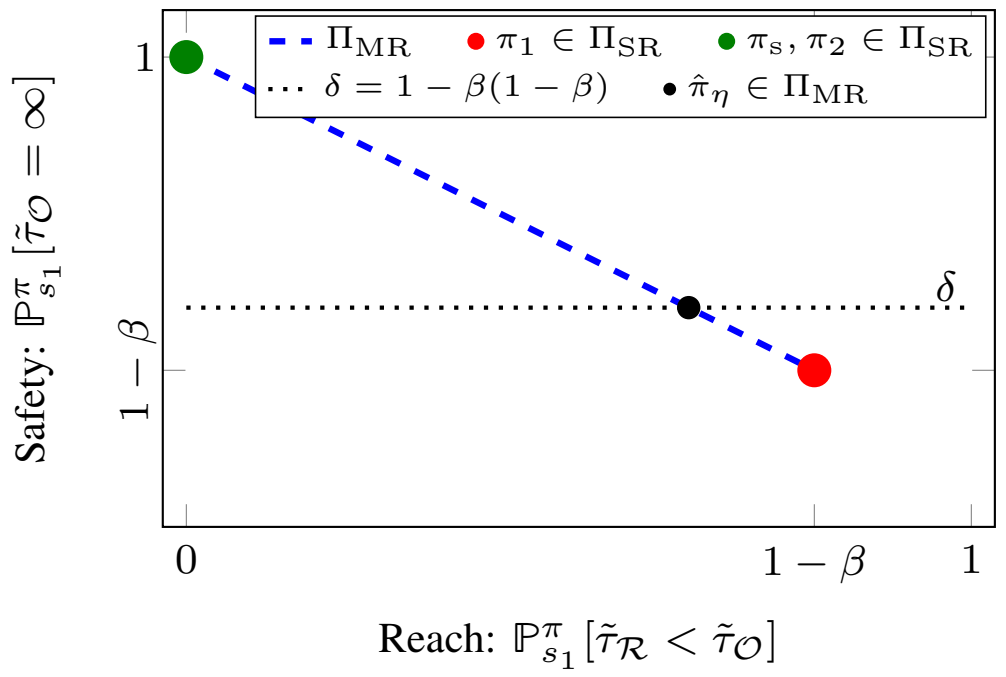
\begin{proof}
\Cref{fig:seg_line} illustrates the construction used in the proof. Given $\varepsilon\in(0,1)$, consider \cref{ex:toy} with $\beta:=1-\sqrt{\varepsilon}\in(0,1)$ and $\delta:=1-\beta(1-\beta)$. We now show that the only feasible $\pi_{\mathrm{s}}\in\PiSR$ is the one with $\pi_{\mathrm{s}}(a_1\mid s_1) = 0$. If $\pi_{\mathrm{s}}(a_1\mid s_1) >0$, then the process eventually reaches $s_3$ with probability $0.1\,\beta \pi_{\mathrm{s}}(a_1\mid s_1)/(0.1\,\pi_{\mathrm{s}}(a_1\mid s_1))=\beta$:
\[
\mathbb P^{\pi_{\mathrm{s}}}_{s_1}[\tilde\tau_{\mathcal O}<\infty]
=\beta>\beta(1-\beta)=1-\delta .
\]
Such a policy is infeasible because it violates the safety constraint. Therefore, the sole feasible $\pi_{\mathrm{s}}\in \PiSR$ satisfies
\[
\mathbb P^{\pi_{\mathrm{s}}}_{s_1}
[\tilde\tau_{\mathcal R}<\tilde\tau_{\mathcal O}]=0,
\]
and is always feasible, since it never visits $s_3$.

Now, define $\pi_1(a_1\mid s_1) = 1$ and $\pi_2(a_1\mid s_1) = 0$ and construct a \emph{lottery} policy $\pi_{\eta} = \pi_1 \oplus_{\eta} \pi_2 \in \PiHR$ that chooses in the first step to play either $\pi_1$ or $\pi_2$ with probabilities $\eta$ and $1-\eta$, respectively. At subsequent visits to $s_1$, the policy $\pi_\eta$ repeats the action chosen in the first step. In addition, there exists a Markov policy $\hat{\pi}_{\eta}\in \PiMR$ with the same state--action marginals at every time as $\pi_{\eta}$, starting from $s_1$~\cite{derman1966note}. Since $\mathcal R$ and $\mathcal O$ are absorbing,
\begin{equation*}
\mathbb{P}^{\pi}_{s_1}[\tilde\tau_{\mathcal X}\le k]
=
\mathbb{P}^{\pi}_{s_1}[\tilde s_k\in\mathcal X]
\end{equation*}
for $\mathcal X\in\{\mathcal R,\mathcal O\}$, every $k\in\Nats$, and any policy $\pi$. Therefore, $\hat{\pi}_{\eta}$ has the same distributions of $\tilde\tau_{\mathcal O}$ and $\tilde\tau_{\mathcal R}$ as $\pi_{\eta}$. Then, choosing $\eta = 1-\beta$:
\begin{align*}
\mathbb{P}^{\hat{\pi}_{\eta}}_{s_1}[\tilde\tau_{\mathcal O}<\infty]
&=(1-\beta)\beta=1-\delta,\\
\mathbb{P}^{\hat{\pi}_{\eta}}_{s_1}[\tilde\tau_{\mathcal R}<\tilde\tau_{\mathcal O}]
&=(1-\beta)^2=\varepsilon .
\end{align*}
Hence, $\hat{\pi}_{\eta} \in \PiMR$ is feasible and outperforms $\pi_{\mathrm{s}}$ by $\varepsilon$.
\end{proof}

\section{Reach--Avoid as a Constrained Average-Reward MDP}
\label{sec:avg}
In this section, we show that the probabilistic reach--avoid problem introduced in~\eqref{pr:RA} can be formulated and solved as a finite constrained multichain average-reward \emph{Markov decision process}~(MDP).

A finite MDP $\left(\mathcal{S},\{\mathcal{A}(s)\}_{s\in\mathcal{S}}, \alpha ,p,r,c\right)$ extends a controlled Markov chain defined in \cref{sec:ra_formulation} where we also have the set of states $\mathcal{S}$, actions $\mathcal{A}$, transition function $p\colon \mathcal{Z} \to \Delta(\mathcal{S})$, and the initial distribution $\alpha\in \Delta(\mathcal{S})$~\cite{Puterman2005,derman1970finite,hernandezlerma1996discrete,feinberg2002handbook}.
In addition, an MDP includes also rewards $r\colon \mathcal{Z} \to \Real$.
\emph{Constrained MDPs} also contain an additional constraint function $c\colon \mathcal{Z} \to \Real$ and a desired threshold $\delta \in \Real$~\cite{altman2021constrained,piunovskiy1997optimal,hordijk1984constrained}.
The cost function $c$ constrains policies that are considered to be feasible solutions.
The objective in an average-reward MDP is known as the \emph{gain} $g_{\mathrm{r}}\colon \mathcal{S} \to \Real$ defined as the following Cesaro limit~\cite[section~8.2]{Puterman2005}:

A finite MDP $\left(\mathcal{S},\{\mathcal{A}(s)\}_{s\in\mathcal{S}}, \alpha,p,r\right)$ extends a controlled Markov chain defined in~\cref{sec:ra_formulation}, with the set of states $\mathcal{S}$, actions $\mathcal{A}$, transition function $p\colon\mathcal{Z}\to\Delta(\mathcal{S})$, and initial distribution $\alpha\in\Delta(\mathcal{S})$~\cite{Puterman2005,derman1970finite,hernandezlerma1996discrete,feinberg2002handbook}. In addition, an MDP includes a reward function $r\colon\mathcal{Z}\to\Real$. \emph{Constrained MDPs} also contain an additional cost function $c\colon\mathcal{Z}\to\Real$~\cite{altman2021constrained,piunovskiy1997optimal,hordijk1984constrained}. The cost function $c$ constrains policies that are considered to be feasible solutions. The objective in an average-reward MDP is known as the \emph{gain}. For a policy $\pi$, the gain $g_{\mathrm r}^{\pi}\colon\mathcal{S}\to\Real$ is defined as~\cite[chapter~8]{Puterman2005}

\begin{align} g_{\mathrm r}^{\pi}(s)
&:=\liminf_{N\to\infty}
\frac{1}{N}\E_s^{\pi}
\left[\sum_{k=0}^{N-1}
r(\tilde{s}_k,\tilde{a}_k)\right].
\label{eq:gain_r}
\end{align}
Similarly, the cost gain
$g_{\mathrm c}^{\pi}\colon\mathcal S\to\Real$
is defined as
\begin{align}
g_{\mathrm{c}}^{\pi}(s)
&\;:=\; 
\liminf_{N\to\infty}
\frac{1}{N} \E_{s}^{\pi} \left[\sum_{k=0}^{N-1} c(\tilde{s}_k, \tilde{a}_k) \right].
\label{eq:gain_c}
\end{align}
In addition, the gains $g_{\mathrm{r}}^{\pi}(\alpha)$ and $g_{\mathrm{c}}^{\pi}(\alpha)$ for the initial distribution $\alpha$ are defined as $ g_f^{\pi}(\alpha) := \sum_{s\in\mathcal{S}}\alpha(s)g_f^{\pi}(s)$ for $f\in\{ \mathrm{r},\mathrm{c}\}$.

With these definitions, the constrained average-reward
problem is formulated as follows:
\begin{equation}\label{pr:RA_gain} \begin{aligned}
    \maximize_{\pi\in\PiHR}\quad & g_{\mathrm{r}}^{\pi}(\alpha) \\
    \stc\quad & g_{\mathrm{c}}^{\pi}(\alpha) \geq \delta.
\end{aligned}
\end{equation}


We now describe how to cast the reach--avoid problem in \cref{def:RA} as a constrained average-reward MDP. The states $\mathcal{S}$, actions $\mathcal{A}$, transition probabilities $p$, and the safety level $\delta$ remain identical to the underlying controlled Markov chain environment and the reach--avoid specification. 
The reward and constraint functions $r,c\colon \mathcal{Z}\rightarrow\Real$ are defined for each $(s,a)\in\mathcal{Z}$ as
\begin{equation}
 \label{eq:reward-cost}
  r(s,a) \;:=\;  \mathbb{1}_{\mathcal{R}}(s), \qquad c(s,a) \;:=\; \mathbb{1}_{\mathcal{S} \setminus \mathcal{O}}(s) = \mathbb{1}_{\mathcal{C}}(s) .
\end{equation}
The following lemma shows the key property of our MDP construction described above. 
\begin{lemma}\label{lemma:gain_prob}
Under~\cref{ass:absorption}, for every $\pi\in\PiHR$ and $s\in\mathcal S$, the gain functions for $r,c$ defined in~\eqref{eq:reward-cost} satisfy
\begin{align}
    g_{\mathrm{r}}^{\pi}(s)
    &=
    \P^{\pi}_{s} \left[\exists k\in\mathbb{N}_0:\, \tilde{s}_k \in \mathcal{R} \right]
    =
    \P^{\pi}_{s} \left[\tilde{\tau}_{\mathcal{R}} < \tilde{\tau}_{\mathcal{O}} \right],
    \label{eq:gain_reach}
    \\
    g_{\mathrm{c}}^{\pi}(s)
    &=
    \P^{\pi}_{s} \left[ \forall k\in\mathbb{N}_0:\, \tilde{s}_k \in \mathcal{C} \right]
    =
    1-\P^{\pi}_{s} \left[\tilde{\tau}_{\mathcal{O}} < \infty \right],
    \label{eq:gain_safe}
\end{align}
and the limits in \eqref{eq:gain_r} and \eqref{eq:gain_c} exist.
\end{lemma}
\begin{proof}
For a fixed policy \(\pi\) and an initial state \(s\), let \(\mathcal L\in\{\mathcal R,\mathcal O\}\). Since $\mathbb{E}$ is a linear operator, for integers $N \ge 1$ we have
\begin{equation*}
\mathbb{E}_s^\pi\!\left[
    \sum_{k=0}^{N-1}\mathbb{1}_{\mathcal{L}}(\tilde{s}_k)
\right]
=
\sum_{k=0}^{N-1}
\mathbb{E}_s^\pi\!\left[
    \mathbb{1}_{\mathcal{L}}(\tilde{s}_k)
\right] \\
=
\sum_{k=0}^{N-1}
\mathbb{P}_s^\pi\!\left[
    \tilde{s}_k \in \mathcal{L}
\right].
\end{equation*}
\Cref{ass:absorption} implies that 
\begin{equation*}
\mathbb P_s^\pi\!\left[\tilde s_k\in\mathcal L\right]
=
\mathbb P_s^\pi\!\left[\tilde\tau_{\mathcal L}\le k\right].
\end{equation*}
As $k\to\infty$, by continuity of probability from below and the preceding equality~(see \citep[theorem~1.1.1(iii)]{durrett2019probability}), we have
\begin{equation*}
  \lim_{k\to \infty} \P^{\pi}_{s}[\tilde{s}_k\in\mathcal{L}] = \mathbb P_s^\pi\!\left[
\bigcup_{k\in\mathbb N_0}
\{\tilde s_k\in\mathcal L\}
\right]=  \P^{\pi}_{s}[\tilde{\tau}_{\mathcal{L}}<\infty].
\end{equation*}
Since a convergent sequence and its Cesàro averages have the same limit, we have
\begin{equation*}
    \lim_{N\rightarrow\infty}
    \frac{1}{N}
    \sum_{k=0}^{N-1}\P^{\pi}_{s}\left[\tilde{s}_k\in\mathcal{L}\right]
    \;=\; 
    \P^{\pi}_{s}\left[\tilde{\tau}_{\mathcal{L}}<\infty\right].
\end{equation*}
If $\mathcal L=\mathcal R$, then $r=\mathbb 1_{\mathcal R}$ and 
\begin{equation*}
g_r^\pi(s)
=
\mathbb P_s^\pi[\tilde\tau_{\mathcal R}<\infty].
\end{equation*}
Due to the fact that $\mathcal O$ is absorbing and $\mathcal R\cap\mathcal O=\varnothing$,
\begin{equation*}
\mathbb P_s^\pi[\tilde\tau_{\mathcal R}<\infty]
=
\mathbb P_s^\pi[\tilde\tau_{\mathcal R}<\tilde\tau_{\mathcal O}].
\end{equation*}
This means whenever $\mathcal R$ is reached at a finite time, it must be reached before $\mathcal O$, almost surely. In terms of $\mathcal C$, we know that
$
\mathbb 1_{\mathcal C}(\tilde s_k)
=
1-\mathbb 1_{\mathcal O}(\tilde s_k).
$ 
Therefore, 
\begin{equation*}
\frac1N\mathbb E_s^\pi\!\left[
\sum_{k=0}^{N-1}\mathbb 1_{\mathcal C}(\tilde s_k)
\right]
=
1-
\frac1N\mathbb E_s^\pi\!\left[
\sum_{k=0}^{N-1}\mathbb 1_{\mathcal O}(\tilde s_k)
\right].
\end{equation*}
Taking $N\to\infty$ gives
\begin{equation*}
\lim_{N\to\infty}
\frac1N\mathbb E_s^\pi\!\left[
\sum_{k=0}^{N-1}\mathbb 1_{\mathcal O}(\tilde s_k)
\right]
=
\mathbb P_s^\pi[\tilde\tau_{\mathcal O}<\infty].
\end{equation*}
Equivalently
\begin{equation*}
g_c^\pi(s)
=
1-\mathbb P_s^\pi[\tilde\tau_{\mathcal O}<\infty].
\end{equation*}
Finally, visiting \(\mathcal O\) means
\begin{equation*}
\mathbb P_s^\pi[
\tilde s_k\in\mathcal C,\ \forall k\in\mathbb N_0]
=
\mathbb P_s^\pi[\tilde\tau_{\mathcal O}=\infty]
=
1-\mathbb P_s^\pi[\tilde\tau_{\mathcal O}<\infty].
\end{equation*}
This completes the proof. 
\end{proof}

With the derivations above, we are now ready to formulate the constrained MDP optimization problem. 
\Cref{lemma:gain_prob} implies the following equivalence theorem.
\begin{theorem}
\label{thm:constrained-optimality}
Under~\Cref{ass:absorption}, let $\Pi_\delta$ be the feasible policy set of~\Cref{pr:RA} in~\Cref{def:RA}. Then
\[
\Pi_\delta
=
\left\{\pi\in\PiHR:
g_{\mathrm c}^\pi(\alpha)\geq\delta\right\},
\]
and, for every $\pi\in\PiHR$,
\[
g_{\mathrm r}^\pi(\alpha)
=
\P_\alpha^\pi
\left[\tilde\tau_{\mathcal R}<\tilde\tau_{\mathcal O}\right].
\]
Consequently, the sets of optimal policies for~\eqref{pr:RA}
and~\eqref{pr:RA_gain} coincide:
\[
\sup_{\pi\in\Pi_\delta}g_{\mathrm r}^\pi(\alpha)
=
\sup_{\pi\in\Pi_\delta}
\P_\alpha^\pi
\left[\tilde\tau_{\mathcal R}<\tilde\tau_{\mathcal O}\right],
\]
and
\[
\operatorname*{arg\,max}_{\pi\in\Pi_\delta}
g_{\mathrm r}^\pi(\alpha)
=
\operatorname*{arg\,max}_{\pi\in\Pi_\delta}
\P_\alpha^\pi
\left[\tilde\tau_{\mathcal R}<\tilde\tau_{\mathcal O}\right].
\]
\end{theorem}

\begin{proof}
By~\Cref{lemma:gain_prob}, averaging over $\alpha$ gives,
for every $\pi\in\PiHR$,
\begin{align*}
g_{\mathrm r}^\pi(\alpha)
&=
\P_\alpha^\pi
\left[\tilde\tau_{\mathcal R}<\tilde\tau_{\mathcal O}\right],\\
g_{\mathrm c}^\pi(\alpha)
&=
1-\P_\alpha^\pi
\left[\tilde\tau_{\mathcal O}<\infty\right].
\end{align*}
These identities establish equality of the feasible
policy sets and pointwise equality of the objectives,
and hence equality of the optimal values and
optimal policy sets.
\end{proof}

By~\cref{thm:constrained-optimality}, the problem reduces to
finding policies $\pi\in\PiHR$ that maximize
$g_{\mathrm{r}}^{\pi}(\alpha)$ while satisfying
$g_{\mathrm{c}}^{\pi}(\alpha)\geq\delta$.
The following section presents the linear programming
formulation of~\eqref{pr:RA_gain} and outlines the procedure
for constructing an optimal policy when the problem is feasible.
\section{Constructing Nonstationary Optimal Policy}
\label{sec:CLP} 
This section presents the linear programming formulation of
Problem~\eqref{pr:RA_gain} and the construction of an optimal
policy when the problem is feasible.

We define the decision variables $\bm{z},\bm{y}\in\Real^{|\mathcal{Z}|}$. For a policy $\pi$, the variable $\bm{z}^{\pi}$ stands for any limit point of the expected long-run state--action frequencies under policy $\pi$. For some convergent subsequence $\{N_j\}$,
\begin{equation}\label{eq:z}
    z^{\pi}_{s,a}
    =
    \lim_{j\rightarrow\infty}
    \frac{1}{N_j}
    \sum_{k=0}^{N_j-1}
    \P^{\pi}_{\alpha}
    \left[
        \tilde{s}_k=s,\ \tilde{a}_k=a
    \right].
\end{equation}
Within the Hordijk--Kallenberg linear
program~\cite{hordijk1984constrained}, the vector $\bm y$
serves as a nonnegative auxiliary variable that links the
initial distribution $\alpha$ to the long-run state--action
frequencies $\bm z$ through the flow-balance
constraints~\cite[theorem~9.4]{kallenberg2011markov}.
Intuitively, it describes how the initial probability mass
is routed through state--action pairs toward the long-run
behavior encoded by $\bm z$. In the lottery-based switching
construction, normalizing its entries over actions at states
with positive auxiliary flow determines the policy used
before switching. Since $\bm y$ represents unnormalized
flow, its components need not sum to one.

Based on \Cref{lemma:gain_prob} and~\eqref{eq:z}, the gains in Problem~\eqref{pr:RA_gain} are linear with respect to $\bm{z}^{\pi}$. In particular,
$g_{\mathrm{r}}^{\pi}(\alpha)=\bm{r}\tr\bm{z}^{\pi}$ and
$g_{\mathrm{c}}^{\pi}(\alpha)=\bm{c}\tr\bm{z}^{\pi}$, where
$\bm{r},\bm{c}\in\Real^{|\mathcal{Z}|}$ are vectors that represent the stage
functions~\eqref{eq:reward-cost} over $\mathcal{Z}$.
The linear program associated with Problem~\eqref{pr:RA_gain} is given below.
\begin{equation} \label{eq:lp_cons}
\begin{aligned}
    \operatorname*{maximize}_{\bm{z}\geq \bm{0},\,
    \bm{y}\geq \bm{0}}
    \quad &
    \bm{r}\tr\bm{z} \\
  \stc
  \quad &
    (\bm{E}-\bm{P})\tr\bm{z}
    =
    \bm{0},\\
    &
    \bm{E}\tr\bm{z}
    +
    (\bm{E}-\bm{P})\tr\bm{y}
    =
    \bm{\alpha},\\
    &
    \bm{c}\tr\bm{z}
    \geq
    \delta.
\end{aligned}
\end{equation}
where $\bm{P}\in\Real^{|\mathcal{Z}|\times|\mathcal{S}|}$
denotes the transition matrix whose row indexed by $(s,a)$
is $p(\cdot\mid s,a)$, and
$\bm{E}\in\Real^{|\mathcal{Z}|\times|\mathcal{S}|}$
is the incidence matrix
\begin{equation}\label{eq:incidence}
    E_{(s,a),j}
    =
    \begin{cases}
        1, & j=s,\\
        0, & j\neq s
    \end{cases}.
\end{equation}
The first two constraints, together with nonnegativity,
characterize the limit points of expected state--action
frequencies achievable by admissible policies, as established
in prior work~\cite{hordijk1984constrained,kallenberg1983linear,kallenberg1994survey}.
Altman and Shwartz~\cite{altman1991markov} also study
state--action frequencies and their role in constrained MDPs.
Projecting this feasible set onto the $\bm z$ variables
yields a bounded polytope.
The third constraint further restricts this polytope
to frequencies that meet the required safety level.
This linear program applies the constrained multichain
formulations of~\cite{hordijk1984constrained,krass2002achieving}
to the reach--avoid problem defined in~\eqref{pr:RA}.

\begin{proposition}
\label{thm:lp_equivalence}
Under~\Cref{ass:absorption}, we have
\begin{equation*}
\exists\,\hat\pi\in\PiHR
\quad\text{such that}\quad
\forall\,\pi\in\PiHR,\qquad
g_{\mathrm c}^{\pi}(\alpha)
\le g_{\mathrm c}^{\hat\pi}(\alpha).
\end{equation*}

Let $\Omega_{\alpha,\delta}$ be the feasible set
of LP~\eqref{eq:lp_cons}. Then
\begin{equation*}
   \Omega_{\alpha,\delta}\neq\varnothing
\iff
g_{\mathrm c}^{\hat\pi}(\alpha)\ge\delta. 
\end{equation*}
Moreover 
\begin{equation*}
\Omega_{\alpha,\delta}\neq\varnothing
\;\Longrightarrow\;
\exists\,(\bm z^\star,\bm y^\star)
\in
\operatorname*{arg\,max}_{(\bm z,\bm y)\in\Omega_{\alpha,\delta}}
\bm r^\top\bm z.
\end{equation*}
For such an optimal solution, there exists a policy
$\pi^\star\in\Pi_\delta$ such that
\begin{equation*}
\begin{aligned}
g_{\mathrm c}^{\pi^\star}(\alpha)
&=\bm c^\top\bm z^\star\ge\delta,\\
g_{\mathrm r}^{\pi^\star}(\alpha)
&=\bm r^\top\bm z^\star
=\max_{\pi\in\Pi_\delta}g_{\mathrm r}^{\pi}(\alpha)
=\max_{\pi\in\Pi_\delta}
\mathbb P_\alpha^\pi[
\tilde\tau_{\mathcal R}<\tilde\tau_{\mathcal O}].
\end{aligned}
\end{equation*}
\end{proposition}

\begin{proof}
For a fix policy $\pi\in\PiHR$ let $\bm z_N^\pi$
presents the expected average state--action frequency
vector over the first $N$ steps. 
Moreover, we define the set of state–action frequency vectors according to~\cite{kallenberg1983linear} as 
\begin{equation*}
\mathcal Q_\alpha
\!:=\!
\left\{
\bm z\ge0\!:\!
\begin{array}{l}
(\bm E-\bm P)^\top\bm z=0,\\[2pt]
\exists\,\bm y\ge0:\quad
\bm E^\top\bm z+(\bm E-\bm P)^\top\bm y=\alpha
\end{array}
\right\}.
\end{equation*}
Theorem~4.7.2 in~\cite{kallenberg1983linear} states
that $\bm z\in\mathcal Q_\alpha$ if and only if
there exists a policy whose expected average
state--action frequency vectors converge to $\bm z$
along a subsequence. That is,
\begin{equation*}
\bm z\in\mathcal Q_\alpha
\iff
\exists\,\pi\in\PiHR,\ \exists\,N_j\uparrow\infty:
\quad
\bm z_{N_j}^{\pi}\longrightarrow\bm z.
\end{equation*}
This convergence, together with~\Cref{lemma:gain_prob},
implies that
\begin{equation*}
\begin{aligned}
g_{\mathrm r}^{\pi}(\alpha)
&=\lim_{j\to\infty}\bm r^\top\bm z_{N_j}^{\pi}
=\bm r^\top\bm z,\\
g_{\mathrm c}^{\pi}(\alpha)
&=\lim_{j\to\infty}\bm c^\top\bm z_{N_j}^{\pi}
=\bm c^\top\bm z.
\end{aligned}
\end{equation*}
The following property 
$
\bm E\bm1=\bm P\bm1=\bm1
$
shows that 
\begin{equation*}
\begin{aligned}
1=\bm1^\top\alpha
&=\bm1^\top\bm E^\top\bm z
+\bm1^\top(\bm E-\bm P)^\top\bm y
=\bm1^\top\bm z.
\end{aligned}
\end{equation*}
The term involving $\bm y$ vanishes because
$(\bm E-\bm P)\bm1=0$.
Together with $\bm z\ge0$, the identity
$\bm1^\top\bm z=1$ shows that $\bm z$ is
a probability vector. Therefore, $\mathcal Q_\alpha$ lies inside the probability simplex. 
We have already shown that this limit belongs to $\mathcal Q_\alpha$. Therefore,
$
\mathcal Q_\alpha\neq\varnothing.
$
Since every component of each vector
$\bm z\in\mathcal Q_\alpha$ lies between $0$ and $1$,
the set $\mathcal Q_\alpha$ is bounded.
As the projection of a polyhedron, $\mathcal Q_\alpha$
is itself a polyhedron and is therefore closed. Thus $\mathcal Q_\alpha$ is nonempty and compact. 

The non-emptiness and compactness of $\mathcal Q_\alpha$,
together with the linearity of the function
$\bm z\mapsto\bm c^\top\bm z$, imply that this function
attains its maximum at some
$\widehat{\bm z}\in\mathcal Q_\alpha$.
By Theorem~4.7.2 in~\cite{kallenberg1983linear},
there exists a policy $\hat\pi\in\PiHR$ whose expected
average frequency vectors have $\widehat{\bm z}$
as a limit point. Hence, by~\Cref{lemma:gain_prob},
\begin{equation*}
g_{\mathrm c}^{\hat\pi}(\alpha)
=
\bm c^\top\widehat{\bm z}
=
\max_{\bm z\in\mathcal Q_\alpha}\bm c^\top\bm z.
\end{equation*}
For any policy $\pi\in\PiHR$, choose a limit point
$\bm z\in\mathcal Q_\alpha$ of $(\bm z_N^\pi)$.
By~\Cref{lemma:gain_prob} and the maximizing property
of $\widehat{\bm z}$,
\begin{equation*}
g_{\mathrm c}^{\pi}(\alpha)
=
\bm c^\top\bm z
\le
\bm c^\top\widehat{\bm z}
=
g_{\mathrm c}^{\hat\pi}(\alpha).
\end{equation*}

Suppose that $\Omega_{\alpha,\delta}\neq\varnothing$.
Then there exists a feasible pair
$(\bm z,\bm y)\in\Omega_{\alpha,\delta}$.
Therefore,
$
\bm z\in\mathcal Q_\alpha$
and $
\bm c^\top\bm z\ge\delta,
$
which implies that 
$
g_{\mathrm c}^{\hat\pi}(\alpha)
=
\bm c^\top\widehat{\bm z}
\ge
\bm c^\top\bm z
\ge\delta.
$ 
On the other side, let us consider
\(g_{\mathrm c}^{\hat\pi}(\alpha)\ge\delta\).
Since \(\widehat{\bm z}\in\mathcal Q_\alpha\),
there exists \(\widehat{\bm y}\ge0\) such that
\((\widehat{\bm z},\widehat{\bm y})\) satisfies
the constraints in~\eqref{eq:lp_cons}.
Together with the preceding safety inequality, this gives
\((\widehat{\bm z},\widehat{\bm y})
\in\Omega_{\alpha,\delta}\).
Therefore, \(\Omega_{\alpha,\delta}\neq\varnothing\).

We define another set that contains the 
$\bm z$ vectors as
\begin{equation*}
\mathcal Q_{\alpha,\delta}
:=
\left\{
\bm z:
\exists\,\bm y\ge0
\text{ such that }
(\bm z,\bm y)\in\Omega_{\alpha,\delta}
\right\}.
\end{equation*}
This set is also nonempty and compact. 
Therefore, there exists
\begin{equation*}
\bm z^\star\in
\operatorname*{arg\,max}_{\bm z\in\mathcal Q_{\alpha,\delta}}
\bm r^\top\bm z.
\end{equation*}
By Theorem~4.7.2 in~\cite{kallenberg1983linear}
and~\Cref{lemma:gain_prob}, we have
\begin{equation*}
\begin{aligned}
g_{\mathrm c}^{\pi^\star}(\alpha)
&=\bm c^\top\bm z^\star\ge\delta,\\
g_{\mathrm r}^{\pi^\star}(\alpha)
&=\bm r^\top\bm z^\star.
\end{aligned}
\end{equation*}
Thus, $\pi^\star\in\Pi_\delta$.
For any $\pi\in\Pi_\delta$, choose a frequency limit
point $\bm z\in\mathcal Q_{\alpha,\delta}$.
By the maximizing property of $\bm z^\star$, we have
\begin{equation*}
g_{\mathrm r}^{\pi}(\alpha)
=
\bm r^\top\bm z
\le
\bm r^\top\bm z^\star
=
g_{\mathrm r}^{\pi^\star}(\alpha).
\end{equation*}
Equivalently 
\begin{equation*}
g_{\mathrm r}^{\pi^\star}(\alpha)
=\max_{\pi\in\Pi_\delta}g_{\mathrm r}^{\pi}(\alpha)
=\max_{\pi\in\Pi_\delta}
\P_\alpha^\pi
[\tilde\tau_{\mathcal R}<\tilde\tau_{\mathcal O}].
\end{equation*}
This completes the proof. 
\end{proof}

Given an optimal LP solution $(\bm z\opt,\bm y\opt)$,
define
\begin{equation*}
z_s\opt := \sum_{a\in\mathcal{A}(s)} z_{s,a}\opt,
\qquad
y_s\opt := \sum_{a\in\mathcal{A}(s)} y_{s,a}\opt.
\end{equation*}
Given an optimal basic solution $(\bm z^*,\bm y^*)$
of~\eqref{eq:lp_cons}, define
$z_s^*:=\sum_{a\in\mathcal A(s)}z_{s,a}^*$ and
$y_s^*:=\sum_{a\in\mathcal A(s)}y_{s,a}^*$.
Consider the stationary policy $\pi_{\mathrm c}$ defined below.

Given an optimal basic solution $(\bm{z}^{*},\bm{y}^{*})$
of~\eqref{eq:lp_cons}, define
$z^{*}_s:=\sum_{a\in\mathcal{A}(s)}z^{*}_{s,a}$ and
$y^{*}_s:=\sum_{a\in\mathcal{A}(s)}y^{*}_{s,a}$.
Consider the case in which the last inequality constraint
$\bm{c}\tr\bm{z}^{*}\geq\delta$ in~\eqref{eq:lp_cons}
is inactive, that is, $\bm{c}\tr\bm{z}^{*}>\delta$.
The solution remains basic after removing this inactive
constraint. It also remains optimal, since a sufficiently
small convex combination with any better unconstrained
solution would remain feasible for~\eqref{eq:lp_cons}
and improve its objective.
The standard multichain average-reward policy construction
in~\cite{avila2021reachability}
and~\cite[theorem~9.40]{kallenberg2016markov}
therefore yields a stationary policy that is optimal
for the unconstrained problem from the initial distribution
$\alpha$, even when $\alpha$ has zero entries.
Since the avoid set is absorbing, this construction
does not decrease the safety gain and therefore also
yields an optimal policy for~\eqref{pr:RA_gain}:
\begin{equation}\label{eq:policy-construction}
    \pi_{\mathrm{c}}(a|s)
    :=
    \begin{cases}
        \dfrac{z^{*}_{s,a}}{z^{*}_s},
        & z^{*}_s>0,\\[2ex]
        \dfrac{y^{*}_{s,a}}{y^{*}_s},
        & z^{*}_s=0,\ y^{*}_s>0,\\[2ex]
        \text{arbitrary},
        & \text{otherwise}.
    \end{cases}
\end{equation}
The stationary construction in~\eqref{eq:policy-construction}
can be applied to any optimal solution of~\eqref{eq:lp_cons};
the optimality guarantee above assumes that the safety
constraint is inactive.
The example used in the proof of~\Cref{thm:stationary_insufficient}
demonstrates that all feasible stationary randomized policies
are strictly suboptimal. Since $\pi_{\mathrm c}\in\Pi_{\mathrm{SR}}$,
this policy cannot be optimal in the scenario considered.
Furthermore, as established by~\Cref{thm:constrained-optimality},
the same conclusion holds for Problem~\eqref{pr:RA_gain}.
Therefore, the stationary construction described
in~\eqref{eq:policy-construction} does not generally ensure optimality.

Moreover, \cref{fig:seg_line} illustrates convex
combinations of the gain pairs $(0,1)$ and $(1-\beta,1-\beta)$,
ordered as $(g_{\mathrm r},g_{\mathrm c})$,
corresponding to the two deterministic stationary policies,
shown for $\beta=0.2$.
For $\delta\in[1-\beta,1]$, the optimal linear programming
gain pair $\bigl((1-\delta)(1-\beta)/\beta,\ \delta\bigr)$
lies on the dashed segment; in particular, the policy
$\hat{\pi}_{\lambda}$ in the proof of
\cref{thm:stationary_insufficient} is optimal at
$\delta=1-\beta(1-\beta)$.
These combinations correspond to convex combinations
of the respective state--action frequency
vectors~\cite{hordijk1984constrained}.

To construct an optimal policy from
$(\bm z\opt,\bm y\opt)$, the target-frequency method
described in~\cite{krass2002achieving} is employed.
This method uses two stationary policies and a single
random switch.

We introduce a second stationary policy $\pi_{\mathrm{a}}$,
used before switching to $\pi_{\mathrm{c}}$, as follows:
\begin{equation}\label{eq:front_policy}
    \pi_{\mathrm{a}}(a|s)
    :=
    \begin{cases}
        \dfrac{y^{*}_{s,a}}{y^{*}_s}, & y^{*}_s>0,\\[2ex]
        \dfrac{z^{*}_{s,a}}{z^{*}_s}, & y^{*}_s=0,\ z^{*}_s>0,\\[2ex]
        \text{arbitrary}, & \text{otherwise}.
    \end{cases}
\end{equation}
In~\eqref{eq:front_policy} and~\eqref{eq:policy-construction},
``arbitrary'' means any fixed probability distribution
on $\mathcal{A}(s)$.

\begin{theorem}\label{thm:opt_policy}
For
$(\bm z\opt,\bm y\opt)\in\Omega_{\alpha,\delta}$,
the lottery policy $f\opt$ defined below induces
$\pi\opt\in\PiHR$ with the same state--action history law
under $\alpha$.
Then $\pi\opt$ is optimal for Problem~\eqref{pr:RA_gain}
and, by~\Cref{thm:constrained-optimality}, for
problem~\eqref{pr:RA}. The lottery policy is
\begin{equation}\label{eq:opt_policy_krass}
    f\opt=\bigl(l,\tilde f\bigr),
\end{equation}
where $l$ is the marker-update rule and $\tilde f$
is the action-selection rule, defined by
\begin{equation}\label{eq:switching_function}
    \tilde f(s,m)
    =
    \begin{cases}
        \pi_{\mathrm{a}}(\cdot\,|s), & m=0,\\[1ex]
        \pi_{\mathrm{c}}(\cdot\,|s), & m=1,
    \end{cases}
\end{equation}
and
\begin{equation}\label{eq:switching-rule}
    l(s,m)
    \!:=\!
    \begin{cases}
        (0,1), & m\!=\!1,\\[1ex]
        (1,0), & m\!=\!0 \land z_s\opt\!=\!0,\\[1ex]
        \left(
            \dfrac{y_s\opt}{z_s\opt+y_s\opt},\!\ \!
            \dfrac{z_s\opt}{z_s\opt+y_s\opt}
        \right),
        & m\!=\!0 \land z_s\opt\!>\!0.
    \end{cases}
\end{equation}
The binary marker is initialized at $m=0$.
At each step, sample its new value from $l(s,m)$
in~\eqref{eq:switching-rule}, with components ordered
as $0,1$, then sample an action from $\tilde f(s,m)$
using the updated marker.
Define $g_q^{f\opt}(\alpha):=g_q^{\pi\opt}(\alpha)$
for $q\in\{\mathrm r,\mathrm c\}$.
\end{theorem}

\begin{proof}
By feasibility of $(\bm z\opt,\bm y\opt)$,
\begin{equation*}
\begin{aligned}
(\bm E-\bm P)\tr\bm z\opt&=\bm 0,\\
\bm E\tr\bm z\opt
+(\bm E-\bm P)\tr\bm y\opt&=\bm\alpha,
\end{aligned}
\end{equation*}
with $\bm z\opt,\bm y\opt\geq\bm 0$.
Thus, the pair satisfies the target-frequency conditions
of~\cite{krass2002achieving}.

The policies~\eqref{eq:policy-construction}
and~\eqref{eq:front_policy}, together with
the switching rule~\eqref{eq:switching-rule},
constitute the lottery-based construction
of Krass and Vrieze~\cite{krass2002achieving}.
This construction attains the target state--action
frequency vector $\bm z\opt$. Therefore, for every
$(s,a)\in\mathcal Z$,
\begin{equation*}
\lim_{N\to\infty}\frac{1}{N}
\sum_{k=0}^{N-1}
\mathbb P_\alpha^{f\opt}
[\tilde s_k=s,\tilde a_k=a]
=z_{s,a}\opt.
\end{equation*}
Moreover, conditioning on the state--action history
gives an equivalent policy $\pi\opt\in\PiHR$
with the same state--action history law under $\alpha$.
By~\cite{derman1966note}, there also exists a Markov
policy $\widehat\pi\opt\in\PiMR$ such that
\begin{equation*}
\mathbb P_\alpha^{\widehat\pi\opt}
[\tilde s_k=s,\tilde a_k=a]
=
\mathbb P_\alpha^{f\opt}
[\tilde s_k=s,\tilde a_k=a]
\end{equation*}
for every $k\in\Nats$ and $(s,a)\in\mathcal Z$.
Consequently, $\pi\opt$ and $\widehat\pi\opt$ attain
the same target frequencies and gains as $f\opt$.
In particular,
\begin{equation*}
\begin{aligned}
g_{\mathrm r}^{\pi\opt}(\alpha)
&=\bm r\tr\bm z\opt,\\
g_{\mathrm c}^{\pi\opt}(\alpha)
&=\bm c\tr\bm z\opt\geq\delta.
\end{aligned}
\end{equation*}
Hence, $\pi\opt$ is feasible for
Problem~\eqref{pr:RA_gain}.
By~\Cref{thm:lp_equivalence},
$\bm r\tr\bm z\opt$ is the optimal value of
Problem~\eqref{pr:RA_gain}, so $\pi\opt$ is optimal.
Finally,~\Cref{thm:constrained-optimality}
implies that $\pi\opt$ is also optimal for
problem~\eqref{pr:RA}.
The same conclusions hold for $\widehat\pi\opt$,
since its gains equal those of $\pi\opt$.
\end{proof}

The key principle behind $f\opt$ is the use of two
stationary policies: $\pi_{\mathrm a}$ governs the initial
phase, while $\pi_{\mathrm c}$ governs the long-run behavior
after switching. The switch occurs at a finite time
almost surely~\cite[lemma~3.1]{krass2002achieving}.

In \cref{ex:toy}, the Pareto frontier consists of a single
line segment connecting the gain pairs of two deterministic
stationary policies. More generally, the frontier may consist
of several such segments. We now illustrate this structure
using a six-state probabilistic reach--avoid example with
three actions.

\begin{example}\label{example_2}
Let $\mathcal S=\{s_1,\ldots,s_6\}$ and
$\mathcal A(s)=\{a_1,a_2,a_3\}$ for every
$s\in\mathcal S$.
The transition matrices corresponding to the three
actions are $P_1=\bm I_{6\times6}$,
\(
P_2=
\begin{bmatrix}
\bm B_1 & \bm B_2\\
\bm0_{3\times3} & \bm B_3
\end{bmatrix},
\)
and
\(
P_3=
\begin{bmatrix}
\bm Q_1 & \bm Q_2\\
\bm0_{3\times3} & \bm Q_3
\end{bmatrix},
\)
where
\begingroup
\small
\setlength{\arraycolsep}{2pt}
\renewcommand{\arraystretch}{0.9}
\begin{equation*}
\begin{aligned}
\bm B_1&=
\begin{bmatrix}
0 & 1 & 0\\
0 & 0 & 0\\
0.17 & 0.17 & 0
\end{bmatrix},
&
\bm Q_1&=
\begin{bmatrix}
0 & 0 & 1\\
0.15 & 0.075 & 0.275\\
0.306 & 0.17 & 0.136
\end{bmatrix},\\[0.5ex]
\bm B_2&=
\begin{bmatrix}
0 & 0 & 0\\
0.64 & 0.35 & 0.01\\
0.34 & 0.12 & 0.20
\end{bmatrix},
&
\bm Q_2&=
\begin{bmatrix}
0 & 0 & 0\\
0 & 0.45 & 0.05\\
0.068 & 0.06 & 0.26
\end{bmatrix},\\[0.5ex]
\bm B_3&=
\begin{bmatrix}
0.60 & 0.10 & 0.30\\
0 & 1 & 0\\
0 & 0 & 1
\end{bmatrix},
&
\bm Q_3&=
\begin{bmatrix}
0.36 & 0.04 & 0.60\\
0 & 1 & 0\\
0 & 0 & 1
\end{bmatrix}.
\end{aligned}
\end{equation*}
\endgroup

The reach and avoid sets are $\mathcal R=\{s_5\}$
and $\mathcal O=\{s_6\}$, respectively.
Since $s_5$ and $s_6$ are absorbing under every action,
\Cref{ass:absorption} holds.
The constraint set is
$\mathcal C=\mathcal S\setminus\mathcal O
=\{s_1,s_2,s_3,s_4,s_5\}$.
Also, we define
\(
r(s,a):=\mathbf{1}_{\mathcal R}(s),\,
c(s,a):=\mathbf{1}_{\mathcal C}(s).
\)
The state reward, cost, and initial distribution vectors
are given by
\(
\bm r_{\mathcal S}=
\begin{bmatrix}
0 \!&\! 0 \!&\! 0 \!&\! 0 \!&\! 1 \!&\! 0
\end{bmatrix}\tr,
\)
\(
\bm c_{\mathcal S}=
\begin{bmatrix}
1 \!&\! 1 \!&\! 1 \!&\! 1 \!&\! 1 \!&\! 0
\end{bmatrix}\tr,
\)
and
\(
\bm\alpha=
\begin{bmatrix}
0.33 \!&\! 0.33 \!&\! 0.34 \!&\! 0 \!&\! 0 \!&\! 0
\end{bmatrix}\tr.
\)

\begin{figure}
    \centering
    \includegraphics[scale=0.35]{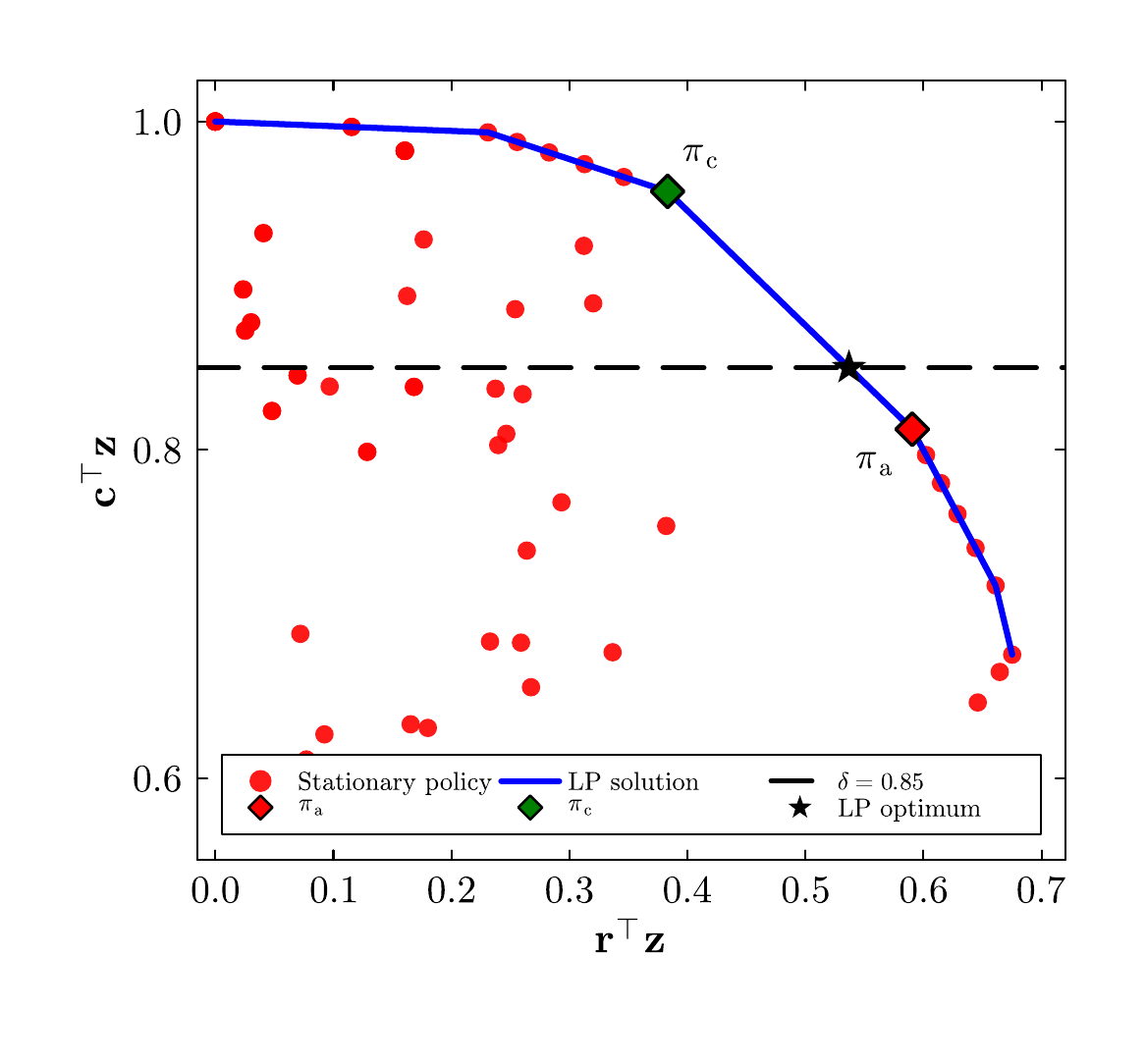}
    \vspace{-0.65cm}
        \caption{Pareto frontier (blue) and stationary policy gain
pairs (red) for the six-state reach--avoid example.
At $\delta=0.85$ (dashed line), the black star marks
the LP optimum attained by $f\opt$.
The green and red diamonds show the gain pairs of
$\pi_{\mathrm c}$ and $\pi_{\mathrm a}$, respectively.}
    \label{fig:pareto_frontier}
\end{figure}

\Cref{fig:pareto_frontier} shows the gain pairs
$(g_{\mathrm r}^\pi(\alpha),g_{\mathrm c}^\pi(\alpha))$.
The red points correspond to deterministic and selected
randomized stationary policies, while the blue curve
shows the Pareto frontier obtained from~\eqref{eq:lp_cons}
as $\delta$ varies.
Some deterministic stationary gain pairs form its vertices,
while other stationary gain pairs lie on its segments or
below it. The segments connecting adjacent vertices
represent convex combinations of the corresponding
state--action frequency vectors~\cite{hordijk1984constrained}.

At $\delta=0.85$, the policies $\pi_{\mathrm c}$ and
$\pi_{\mathrm a}$ constructed from the same optimal LP
solution by~\eqref{eq:policy-construction}
and~\eqref{eq:front_policy}, respectively,
\(
(g_{\mathrm r}^{\pi_{\mathrm c}}(\alpha),
 g_{\mathrm c}^{\pi_{\mathrm c}}(\alpha))
\approx(0.3832,0.9574)
\)
and
\(
(g_{\mathrm r}^{\pi_{\mathrm a}}(\alpha),
 g_{\mathrm c}^{\pi_{\mathrm a}}(\alpha))
\approx(0.5905,0.8125)
\).
The switching policy $f\opt$ uses $\pi_{\mathrm a}$
before switching permanently to $\pi_{\mathrm c}$
according to~\eqref{eq:switching-rule}.
Since $f\opt$ attains the target frequency vector
$\bm z\opt$, its gains satisfy
$g_{\mathrm r}^{f\opt}(\alpha)
=\bm r\tr\bm z\opt\approx0.5369$
and
$g_{\mathrm c}^{f\opt}(\alpha)
=\bm c\tr\bm z\opt=0.85$.
Thus, $f\opt$ attains the LP optimum.
The stationary policy $\pi_{\mathrm c}$ is feasible
but has a strictly smaller reach probability, whereas
$\pi_{\mathrm a}$ violates the safety constraint.
\end{example}

\section{Related Work}
\label{sec:rp}
Related work includes probabilistic reach--avoid
synthesis~\cite{wicker2024probabilistic} and safe exploration
in constrained MDPs~\cite{ni2025safe}.
This section focuses on the formulations proposed by
{\'A}vila and Junca~\cite{avila2021reachability}
and Akametalu et al.~\cite{akametalu2023minimum},
because their reward formulations are closely related
to our reach--avoid analysis.
It also presents the connection between the gain function
and probabilistic safety level sets for the avoidance
problem, a special case of~\eqref{pr:RA}.
Furthermore, the minimum discounted reward (MDR)
Hamilton--Jacobi (HJ) formulation and the average-reward
formulation are compared to examine the effect of
the discount factor on safety level sets.

\subsection{Stationary Average-Reward Formulations}

The stationary policy constructed
in~\cite{avila2021reachability}, denoted by
$\pi_{\mathrm c}$, follows the standard multichain
average-reward policy construction
in~\cite{Puterman2005}.
The counterexample in~\Cref{sec:ra_formulation}
shows that $\pi_{\mathrm c}$ can be suboptimal
for Problem~\eqref{pr:RA_gain}.

The counterexample in the proof
of~\Cref{thm:stationary_insufficient} shows that all
feasible stationary policies are strictly suboptimal
for problem~\eqref{pr:RA} in that example.
This includes any feasible stationary policy obtained
from the Lagrangian formulation
of~\cite{avila2021reachability}.
By~\Cref{thm:constrained-optimality}, the same conclusion
holds for Problem~\eqref{pr:RA_gain}.

{\'A}vila and Junca~\cite{avila2021reachability} investigated the reach--avoid problem using an average-reward MDP and applied Lagrangian relaxation to address the avoid set. In their framework, the long-run average reward quantifies the probability of eventually reaching the target set, while a penalty term enforces the avoid-set constraint. The primary distinction from the present approach lies in how the avoid set is modeled. When the avoid set is treated as an absorbing cemetery set, any history that enters it cannot subsequently reach the target set. Hence,
\[
\P_{\alpha}^{\pi}
\left[
\tilde{\tau}_{\mathcal{R}}<\infty
\right]
=
\P_{\alpha}^{\pi}
\left[
\tilde{\tau}_{\mathcal{R}}<\tilde{\tau}_{\mathcal{O}}
\right].
\]
Therefore, under the cemetery-set modeling assumption,
their reach--avoid formulation coincides with
problem~\eqref{pr:RA}.
This work investigates whether the stationary policy
$\pi_{\mathrm c}$ derived from the average-reward
formulation attains the optimal constrained
reach--avoid performance.

\subsection{Cemetery Avoid Sets}

{\'A}vila and Junca~\cite{avila2021reachability} formulated the following optimization problem:
\begin{equation}
\label{pr:avila2}
\begin{aligned}
    \max_{\pi\in\PiHR}\quad &
    \mathbb{P}_{\alpha}^{\pi}
    \left[
        \tilde{\tau}_{\mathcal{R}}<\infty
    \right] \\
    \stc\quad &
    \mathbb{P}_{\alpha}^{\pi}
    \left[
        \tilde{\tau}_{\mathcal{O}}<\infty
    \right]
    \leq 1-\delta .
\end{aligned}
\end{equation}
The primary distinction between \eqref{pr:avila2} and Problem \eqref{pr:RA} concerns the treatment of the avoid set. In \eqref{pr:avila2}, the avoid set $\mathcal{O}$ is not required to be absorbing, which allows a history to enter $\mathcal{O}$ and subsequently reach the target set. In contrast, Problem \eqref{pr:RA} considers a history successful only if it reaches the target set prior to entering $\mathcal{O}$, that is, when $\tilde{\tau}_{\mathcal{R}}<\tilde{\tau}_{\mathcal{O}}$.

If the avoid set in~\eqref{pr:avila2} is treated as an absorbing cemetery set, any history that enters $\mathcal{O}$ can no longer reach $\mathcal{R}$. Under this cemetery-set assumption, the two reach--avoid objectives become equivalent. Consequently, the distinction between~\eqref{pr:avila2} and problem~\eqref{pr:RA} arises primarily from differences in problem modelling, rather than from a fundamental difference in the reach--avoid objective.

According to~\Cref{ass:absorption}, $\mathcal O$ serves
as an absorbing set and functions as a cemetery set
within this framework. The safety constraint in problem~\eqref{pr:RA} is explicitly represented by the average-cost constraint $g_{\mathrm{c}}^{\pi}(\alpha)\geq\delta$. Consequently, Problem~\eqref{pr:RA_gain} can be formulated as a standard constrained multichain average-reward Markov decision process~\cite{hordijk1984constrained,altman2021constrained}. This approach eliminates the need to expand the model state space to verify whether $\mathcal{O}$ has been visited. However, the optimal policy described in~\Cref{thm:opt_policy} incorporates a binary indicator $m\in\{0,1\}$ alongside the state $s$ to indicate whether the policy has transitioned to its long-run mode.

\subsection{Safety Level Sets}

In this section, we consider a special case of problem~\eqref{pr:RA} in which only the safety constraint matters. Instead of maximizing the probability of reaching $\mathcal{R}$ before $\mathcal{O}$, we just maintain the prescribed safety level. This turns the problem into an avoidance problem:
\begin{equation}\label{pr:avoid}
    \exists\,\pi\in\PiHR
    \ \ \text{such that}\ \
    \P^{\pi}_{s}
    \left[
        \tilde{s}_k\in\mathcal{C},\ \forall k\in\mathbb{N}_0
    \right]
    \geq
    \delta.
\end{equation}

Let us begin with the following definition.
\begin{definition}\label{def:safe_set}
For a safety level $\delta\in(0,1]$, the probabilistic safe set is
\begin{equation}\label{eq:safe_set}
    \mathcal{K}_{\delta}
    \!:=\!
    \left\{
        s\in\mathcal{C}
        \,\!:\!\,
        \exists\,\pi\in\PiHR
        \ \text{s.t.}\
        \P^{\pi}_{s}
        \left[
            \tilde{s}_k\in\mathcal{C},\ \forall k\!\in\!\mathbb{N}_0
        \right]
        \!\geq\!\delta
    \right\}.
\end{equation}
\end{definition}

In this section, we are interested in the connection between the
safety gain $g_{\mathrm c}^{\pi}(s)$ and $\mathcal K_\delta$.
Note that for $\delta=1$, the set $\mathcal K_1$ is the
almost-sure control invariant set: the set of initial states
from which the history can be kept within $\mathcal C$
for all time with probability
one~\cite{avila2021reachability,abate2008probabilistic}.

By~\Cref{pro:safe_set}, a single stationary policy maximizes
the safety gain from every initial state.
Hence, by~\Cref{lemma:gain_prob}, condition~\eqref{pr:avoid}
holds if and only if $g_{\mathrm c}^{\pi}(s)\geq\delta$
for some $\pi\in\PiHR$, or equivalently,
\begin{equation}\label{eq:opt_safety_gain}
    g_{\mathrm c}\opt(s)
    :=
    \sup_{\pi\in\PiHR}g_{\mathrm c}^{\pi}(s)
    \geq\delta.
\end{equation}
The condition $g_{\mathrm c}\opt(s)\geq\delta$
concerns a fixed initial state $s$.
For an initial distribution $\alpha$,
by~\Cref{thm:lp_equivalence},
LP~\eqref{eq:lp_cons} is feasible if and only if
\begin{equation*}
\sup_{\pi\in\PiHR}g_{\mathrm c}^{\pi}(\alpha)
=
\sum_{s\in\mathcal S}\alpha(s)g_{\mathrm c}\opt(s)
\geq\delta.
\end{equation*}
Thus, LP feasibility depends on the weighted optimal safety
gain and does not require
$\operatorname{supp}(\alpha)\subseteq\mathcal K_\delta$.

The reach--avoid problem subsequently maximizes the
reach probability over the policies that satisfy the
prescribed safety level $\delta$.
The avoidance problem therefore characterizes the states
from which Problem~\eqref{pr:RA_gain} is feasible.
By~\Cref{lemma:gain_prob}, $g_{\mathrm c}^{\star}(s)=0$
for every $s\in\mathcal O$.
The following proposition states this characterization
formally and establishes attainment of the optimal
safety gain.

\begin{proposition}\label{pro:safe_set}
There exists a stationary policy $\pi\opt\in\PiSR$
that maximizes the safety gain from every initial state.
For every $s\in\mathcal S$ and $\delta\in(0,1]$,
\begin{equation*}
s\in\mathcal K_\delta
\iff
g_{\mathrm c}\opt(s)
=
g_{\mathrm c}^{\pi\opt}(s)
\geq\delta.
\end{equation*}
In particular, $s\in\mathcal K_1$ if and only if
$g_{\mathrm c}\opt(s)=1$.
\end{proposition}

\begin{proof}
By~\cite[theorem~9.1.8, p.~451]{Puterman2005},
there exists a deterministic stationary policy
$\pi\opt\in\PiSR$ satisfying
$g_{\mathrm c}^{\pi\opt}(s)=g_{\mathrm c}\opt(s)$
for every $s\in\mathcal S$.
Since $0\leq c\leq1$, we have
$0\leq g_{\mathrm c}\opt(s)\leq1$.
By~\Cref{def:safe_set} and~\Cref{lemma:gain_prob},
\begin{equation*}
\begin{aligned}
s\in\mathcal K_\delta
&\iff \exists\,\pi\in\PiHR:
       g_{\mathrm c}^{\pi}(s)\geq\delta\\
&\iff g_{\mathrm c}\opt(s)\geq\delta\\
&\iff g_{\mathrm c}^{\pi\opt}(s)\geq\delta.
\end{aligned}
\end{equation*}
The second equivalence uses attainment of the supremum
by $\pi\opt$.
For $\delta=1$, the bound
$g_{\mathrm c}\opt(s)\leq1$ gives the stated
characterization of $\mathcal K_1$.
\end{proof}

To compute $g_{\mathrm{c}}^{\star}$, and thereby $\mathcal{K}_{\delta}$, we
consider the following linear
program~\cite[\S 9.3]{Puterman2005}:
\begin{equation}\label{eq:lp_gain_bias}
\begin{aligned}
    \minimize_{\bm{g}\in\Real^{|\mathcal{S}|},\,
    \bm{h}\in\Real^{|\mathcal{S}|}}
    \quad &
    \bm{e}\tr\bm{g} \\
    \stc  
    \quad &
    (\bm{E}-\bm{P})\,\bm{g}
    \geq
    \bm{0},\\
    &
    \bm{E}\,\bm{g}
    +
    (\bm{E}-\bm{P})\,\bm{h}
    \geq
    \bm{c},
\end{aligned}
\end{equation}
where $\bm{e}\in\Real^{|\mathcal{S}|}$ denotes the all-ones vector.
More generally, any strictly positive weight vector can be used in the
objective, yielding the same componentwise-minimal solution
$\bm{g}=\bm{g}_{\mathrm{c}}^{\star}$~\cite[\S 9.3]{Puterman2005}.
The matrices $\bm{P}$ and $\bm{E}$ are the same matrices as
in~\eqref{eq:lp_cons}, here untransposed; each constraint row
corresponds to one pair $(s,a)\in\mathcal{Z}$.
Any optimal solution $(\bm{g}^{\star},\bm{h}^{\star})$ satisfies
$g^{\star}(s)=g_{\mathrm{c}}^{\star}(s)$ for all
$s\in\mathcal{S}$~\cite[\S 9.3]{Puterman2005}.
By~\Cref{pro:safe_set}, the solution
of~\eqref{eq:lp_gain_bias} yields the safe set
$\mathcal{K}_{\delta}
=\{s\in\mathcal{C}:g^{\star}(s)\geq\delta\}$
for every $\delta\in(0,1]$.

\subsection{Discounted Approximation}
For the special case $\delta=1$, the average-reward
MDP formulation characterizes the almost-sure controlled
invariant set $\mathcal K_1$.
Akametalu et al. proposed a minimum discounted reward
formulation for approximating reachable sets and derived
a time-independent HJ equation for its value
function~\cite{akametalu2023minimum}.
In this work, we compare the minimum discounted reward
approximation with the average-reward characterization
of $\mathcal K_1$ for different discount factors.
For the pure avoidance problem with $\delta=1$,
an optimal stationary policy exists, as stated
in~\Cref{pro:safe_set}.

According to~\Cref{pro:safe_set}, if $s\in\mathcal K_1$,
there exists an optimal stationary policy that keeps
the history in $\mathcal C$ almost surely.
If $s\notin\mathcal K_1$, then
$g_{\mathrm c}^{\star}(s)<1$, so every policy has
a positive probability of eventually leaving
the constraint set.

The discount factor introduced
in~\cite{akametalu2023minimum} is
$\gamma=\exp(-\lambda\Delta t)$, where $\Delta t$
is the time step and $\lambda$ is the discount rate.
Setting $\lambda=0$ gives $\gamma=1$.
In our numerical comparison in~\Cref{sec:nv},
the MDR Bellman iteration is also evaluated with
$\gamma=1$ as the undiscounted case, although the
contraction guarantee associated with $\gamma<1$
no longer applies.

\section{Numerical Validation}\label{sec:nv}
This section presents the numerical validation of the
proposed framework using two illustrative examples.
The avoidance problem described in~\cref{sec:rp}
is demonstrated using a nonlinear inverted pendulum.
The optimal safety gain $g_{\mathrm c}^{\star}$
is computed using the linear program~\eqref{eq:lp_gain_bias},
and its superlevel sets characterize the probabilistic
safe sets $\mathcal K_\delta$ across the state space
(\Cref{pro:safe_set}).
Moreover, the minimum discounted reward approximation
is compared with the average-reward characterization
of $\mathcal K_1$.
Second, the reach--avoid problem~\eqref{pr:RA}
is validated using a stochastic variant of Zermelo's
navigation problem.
The constrained linear program~\eqref{eq:lp_cons}
is solved for a range of safety levels $\delta$.
The lottery-based policy $f^{*}$
of~\eqref{eq:opt_policy_krass} is compared with the
stationary policy $\pi_{\mathrm c}$
of~\eqref{eq:policy-construction}, which is constructed
from the same LP solution.
Both continuous systems are approximated by finite MDPs
as described below.

All experiments are implemented in Julia, utilizing
JuMP and MOSEK to solve the linear programs.\footnote{
Computations for the inverted pendulum are executed
on a workstation equipped with a 32-core AMD
Threadripper 3970X CPU and 256 GB of RAM.
Zermelo's navigation problem is solved on a MacBook Pro
with an M4 Pro chip and 24 GB of RAM.}
Both systems have continuous-state dynamics of the form
$\tilde{s}_{k+1}=F(\tilde{s}_k,\tilde{a}_k,\tilde{d}_k)$,
where $\tilde{d}_k$ denotes a random disturbance.
The systems are approximated by finite-state MDPs
on uniform state and action grids.
For each discrete state--action pair
$(s,a)\in\mathcal Z$, the transition function
is estimated using Monte Carlo simulation.
Specifically, 100 disturbance samples are drawn.
The continuous successor states are computed using
the system dynamics and mapped to their nearest grid
points using a KD-tree nearest-neighbor
search\footnote{The nearest-neighbor search is used
only to map continuous successor states to the
discrete grid; it is not used for value approximation.},
so that $p(s'\mid s,a)$ represents the empirical
fraction of samples mapped to $s'$.
Grid points located within $\mathcal R$ and
$\mathcal O$ are designated as absorbing in accordance
with~\Cref{ass:absorption}.
The safety and optimality guarantees apply to the
estimated finite MDPs and do not automatically extend
to the original continuous-state systems.

\emph{Inverted Pendulum:} The nonlinear pendulum dynamics are
\begin{equation}\label{eq:IP_system}
\begin{aligned}
\theta_{k+1} &= \theta_k + \omega_k\,\Delta t, \\
\omega_{k+1} &= \omega_k
+ \left(
\frac{g_p}{l_p}\sin\theta_k
+ \frac{1}{m_p l_p^{2}}\,u_k
+ \tilde{d}_k
\right)\Delta t,
\end{aligned}
\end{equation}
where the state is $s_k=(\theta_k,\omega_k)$, with $\theta$
the angle and $\omega$ the angular velocity, and $u_k$
is the control torque.
The physical parameters are the gravitational acceleration
$g_p=9.81\,$m/s$^2$, the pendulum length $l_p=1\,$m,
and the mass $m_p=1\,$kg; the time step is
$\Delta t=0.1\,\mathrm{s}$.
The disturbance $\tilde d_k$ is drawn from
$\mathcal N(0,1)$ truncated to $[-0.75,0.75]$,
independently across time steps.
The constraint set is
$\mathcal C=\{(\theta,\omega)\in\mathcal S:
|\theta|\leq0.3,\ |\omega|\leq0.6\}$,
and the control input is restricted to $u_k\in[-3,3]$.
The state space is discretized on a uniform
$201\times201$ grid extending beyond $\mathcal C$.
The grid points outside $\mathcal C$ form the avoid set
$\mathcal O=\mathcal S\setminus\mathcal C$ and are
designated as absorbing.
The input interval is discretized into $81$ equally
spaced actions, yielding
$|\mathcal Z|=201^2\times81\approx3.3\times10^6$
state--action pairs.

\begin{figure}
    \centering
    \includegraphics[height=5cm]{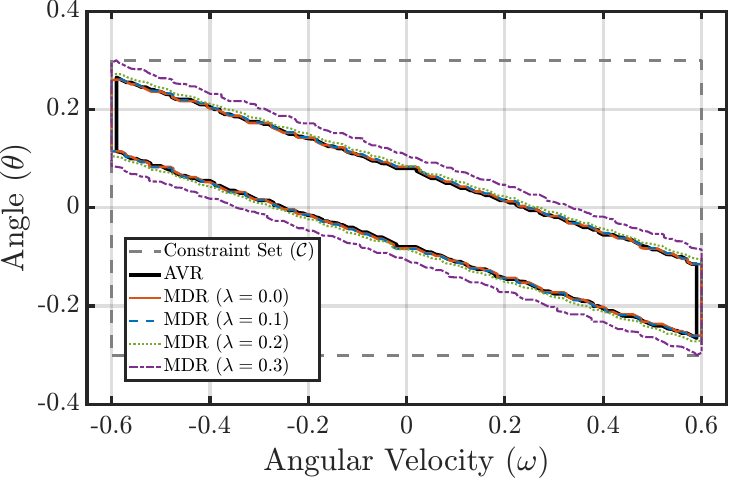}
    \caption{Comparison of $\mathcal K_1$ computed by
    LP~\eqref{eq:lp_gain_bias} with the MDR approximations
    for several discount rates $\lambda$ in the inverted
    pendulum example. The undiscounted case $\lambda=0$
    ($\gamma=1$) gives the closest approximation.}
    \label{fig:comparison}
\end{figure}

\Cref{fig:comparison} compares $\mathcal K_1$ computed
using the average-reward LP~\eqref{eq:lp_gain_bias}
with the MDR approximations for different discount rates.
The closest MDR approximation is obtained for $\lambda=0$,
which corresponds to $\gamma=1$.
For $\lambda=0.1$ and $\lambda=0.2$, the MDR approximation
extends slightly beyond $\mathcal K_1$.
The expansion is more pronounced for $\lambda=0.3$.
States included in the MDR approximation but outside
$\mathcal K_1$ have a positive probability of leaving
the constraint set under every policy in the estimated
finite MDP.
These results illustrate sensitivity to the discount factor
in approximating the solution to Problem~\eqref{pr:avoid}.

\begin{figure}
    \centering
    \includegraphics[height=5cm]{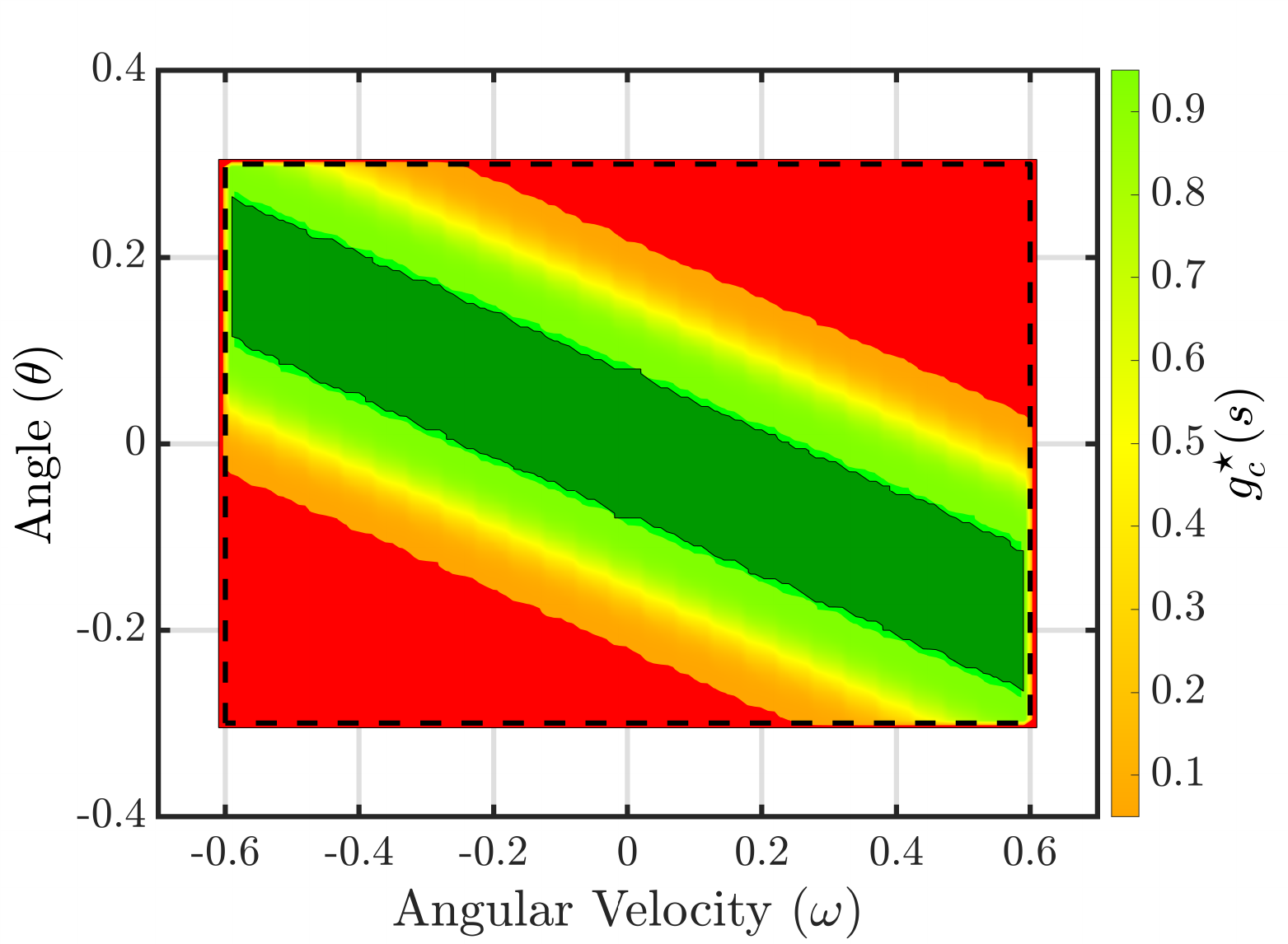}
    \vspace{-.2cm}
    \caption{Contour plot of the optimal safety gain
    $g_{\mathrm{c}}^{\star}(s)$ for the inverted pendulum, computed
    by~\eqref{eq:lp_gain_bias}. Red regions are
    $g_{\mathrm{c}}^{\star}(s)=0$, and the dark green region is the
    almost-sure safe set $\mathcal{K}_{1}$.}
    \label{fig:gain_IP}
\end{figure}

\Cref{fig:gain_IP} illustrates the optimal safety gain
$g_{\mathrm{c}}^{\star}(s)$, which is computed through~\eqref{eq:lp_gain_bias}. The dashed rectangle indicates the constraint set
$\mathcal{C}$. According to
\Cref{pro:safe_set}, the contour levels of
$g_{\mathrm{c}}^{\star}$ show different sets of $\{\mathcal{K}_{\delta}\}_{\delta\in(0,1]}$. States depicted in
red represent $g_{\mathrm{c}}^{\star}(s)=0$. The dark green region represents the
almost-sure safe set $\mathcal{K}_{1}$, where
$g_{\mathrm{c}}^{\star}(s)=1$. The
intermediate orange-to-light-green band comprises states with
$g_{\mathrm{c}}^{\star}(s)\in(0,1)$. According to~\Cref{lemma:gain_prob},
the safety gain at such a state represents the optimal probability of
remaining in $\mathcal{C}$ for all future time.

\emph{Zermelo's Navigation Problem:}
We consider Zermelo's navigation problem in a reach--avoid setting, where a ship moves through a bounded river toward a port under stochastic wind disturbances. The riverbanks establish the navigational boundaries, while the island and regions outside the river represent areas to be avoided. Whenever the ship enters any of these regions, it remains there indefinitely. A shallow-water region is also available, providing a safe location for the ship to anchor. The primary objective is to maximize the probability of reaching the port before entering the avoid regions, while ensuring that the ship remains within the prescribed boundaries with a specified probability.

The state is $s_k=(\mathrm{x}_k,\mathrm{y}_k)$, where
$\mathrm{x}_k$ and $\mathrm{y}_k$ are the position
coordinates at time step $k$. The discrete-time kinematics are
\begin{equation}\label{eq:ki_ship}
\begin{aligned}
\mathrm{x}_{k+1}
&=\mathrm{x}_k+\Delta t
  \left(v\cos u_k+\tilde{w}_k^\mathrm{x}\right),\\
\mathrm{y}_{k+1}
&=\mathrm{y}_k+\Delta t
  \left(v\sin u_k+\tilde{w}_k^\mathrm{y}\right),
\end{aligned}
\end{equation}
where $v=1$ is the constant speed, $u_k$ is the steering
angle in radians, and $\Delta t=0.1$.
The disturbances are independent across coordinates
and time, with
$\tilde{w}_k^\mathrm{x}\sim\mathcal N(-0.5,1)$ and
$\tilde{w}_k^\mathrm{y}\sim\mathcal N(-0.75,1)$;
the second parameter denotes the variance.
The continuous state space $[0,10]\times[0,5]$
is discretized using a uniform grid with step size $0.1$.
The steering input is discretized over $[-1,1]$ radians
with step size $0.05$ radians.

A circular island of radius $0.75$ is centered at
$(5,2.25)$, with
$\mathcal O_{\mathrm{island}}
=\{s\in\mathcal S:\|s-(5,2.25)\|_2\leq0.75\}$.
The river constraint set is
$\mathcal C_{\mathrm{river}}
=\{(\mathrm{x},\mathrm{y})\in\mathcal S:
1<\mathrm{x}<9,\ 1<\mathrm{y}<4\}$.
The constraint and avoid sets are
$\mathcal C
=\mathcal C_{\mathrm{river}}\setminus\mathcal O_{\mathrm{island}}$
and
$\mathcal O
=\mathcal O_{\mathrm{island}}
 \cup(\mathcal S\setminus\mathcal C_{\mathrm{river}})$,
respectively.
The reach set is
$\mathcal R
=\{(\mathrm{x},\mathrm{y})\in\mathcal S:
8\leq\mathrm{x}\leq8.95,\ 1.1\leq\mathrm{y}\leq2.1\}$.

The shallow-water region is an axis-aligned ellipse
centered at $(4,3.25)$, with semi-axes
$a_{\mathrm{sw}}=0.8$ and $b_{\mathrm{sw}}=0.25$ in the
$\mathrm{x}$- and $\mathrm{y}$-directions, respectively.
An additional anchoring action $a_{\mathrm{anc}}$
is feasible only within this region and keeps the ship
at its current position deterministically:
$p(s\mid s,a_{\mathrm{anc}})=1$.
The initial state is $s_0=(2,2.5)$.
\begin{figure}[t]
    \centering
\begin{tikzpicture}
\begin{axis}[
    width=\columnwidth, height=4.8cm,
    xmin=0, xmax=0.7, ymin=0.65, ymax=1,
    xtick={0,0.1,0.2,0.3,0.4,0.5,0.6},
    ytick={0.65,0.70,0.75,0.80,0.85,0.90,0.95,1.00},
    tick label style={font=\footnotesize},
    xticklabel style={/pgf/number format/fixed,
        /pgf/number format/precision=1,
        /pgf/number format/fixed zerofill},
    yticklabel style={/pgf/number format/fixed,
        /pgf/number format/precision=2,
        /pgf/number format/fixed zerofill},
    xlabel={$\mathbf r^\top\mathbf z$},
    ylabel={$\mathbf c^\top\mathbf z$},
    label style={font=\normalsize},
    xlabel style={yshift=2pt},
    ylabel style={yshift=-2pt},
    axis line style={line width=0.4pt},
    tick style={line width=0.35pt}, tick align=inside,
    scaled x ticks=false, scaled y ticks=false,
    grid=major,
    grid style={gray,opacity=0.18,line width=0.3pt},
    legend style={
        at={(0.035,0.055)}, anchor=south west,
        font=\footnotesize, draw=black, fill=white,
        line width=0.35pt, inner sep=2pt,
        row sep=0pt, cells={anchor=west}
    },
    execute at end axis={
        \draw[black!60,densely dotted,line width=0.3pt,
              shorten <=2.5pt]
            (axis cs:0.1379655846,0.9375395826) --
            (axis cs:0.1379655846,0.9532895826);
        \draw[black!60,densely dotted,line width=0.3pt,
              shorten <=2.5pt]
            (axis cs:0.1599808129,0.9322000000) --
            (axis cs:0.1599808129,0.9532895826);
        \draw[black,line width=0.45pt,|-|]
            (axis cs:0.1379655846,0.9532895826) --
            (axis cs:0.1599808129,0.9532895826);
        \node[align=center,fill=white,inner sep=1pt,
              font=\scriptsize] (gaplabel1)
            at (axis description cs:0.40,0.89)
            {$\delta=0.9322$\\$+15.96\%$};
        \draw[black!60,line width=0.3pt]
            (gaplabel1.west) --
            (axis cs:0.1489731988,0.9532895826);
        \draw[black!60,densely dotted,line width=0.3pt,
              shorten <=2.5pt]
            (axis cs:0.5114806894,0.7897626138) --
            (axis cs:0.5114806894,0.8055126138);
        \draw[black!60,densely dotted,line width=0.3pt,
              shorten <=2.5pt]
            (axis cs:0.5552890624,0.7613999994) --
            (axis cs:0.5552890624,0.8055126138);
        \draw[black,line width=0.45pt,|-|]
            (axis cs:0.5114806894,0.8055126138) --
            (axis cs:0.5552890624,0.8055126138);
        \node[align=center,fill=white,inner sep=1pt,
              font=\scriptsize] (gaplabel2)
            at (axis description cs:0.80,0.63)
            {$\delta=0.7614$\\$+8.57\%$};
        \draw[black!60,line width=0.3pt]
            (gaplabel2.south) --
            (axis cs:0.5333848759,0.8055126138);
        %
        \node[above=4pt,font=\scriptsize,fill=white,
              inner sep=0.5pt]
            at (axis cs:0.0828593143,0.9499999906)
            {$\delta=0.95$};
        \node[left=4pt,font=\scriptsize,fill=white,
              inner sep=0.5pt]
            at (axis cs:0.6691070157,0.6691070157)
            {$\delta=0.65$};
    }
]
\addplot[draw={rgb,255:red,0;green,157;blue,255},line width=1.5pt] coordinates {(0.082859323707,0.950000000000) (0.159980812948,0.932200000000) (0.188123427330,0.925000000000) (0.268582990266,0.900000000000) (0.340207734417,0.875000000000) (0.401326558459,0.850000000000) (0.450968167850,0.825000000000) (0.495600777145,0.800000000000) (0.534282774856,0.775000000000) (0.555289062993,0.761400000000) (0.571436286305,0.750000000000) (0.582639115897,0.742000000000) (0.604370132506,0.725000000000) (0.634201245561,0.700000000000) (0.662448769784,0.675000000000) (0.669002195408,0.669200000000) (0.669107015752,0.669107015752)};
\addlegendentry{Optima of~\eqref{eq:lp_cons}}
\addplot[only marks,mark=*,mark size=2.5pt,mark options={draw={rgb,255:red,227;green,26;blue,28},fill={rgb,255:red,227;green,26;blue,28},line width=0.45pt}] coordinates {(0.669107015701,0.669107015701) (0.555289062438,0.761399999445) (0.159980812899,0.932199999952) (0.082859314345,0.949999990637)};
\addlegendentry{$f\opt$}
\addplot[only marks,mark=diamond*,mark size=2.7pt,mark options={draw=black,fill={rgb,255:red,0;green,138;blue,46},line width=0.45pt}] coordinates {(0.669107015701,0.669107015701) (0.511480689443,0.789762613826) (0.137965584639,0.937539582592) (0.082859314345,0.949999990637)};
\addlegendentry{$\pi_{\mathrm c}$}
\end{axis}
\end{tikzpicture}
\caption{Safety--performance trade-off for Zermelo's navigation
problem. The blue line connects sampled optima
of~\eqref{eq:lp_cons}; red circles and green diamonds
show independently evaluated gains of $f\opt$ and
$\pi_{\mathrm c}$, respectively. Percentages show relative
reach probability increases over $\pi_{\mathrm c}$.}
\label{fig:lp_ZN}
\end{figure}
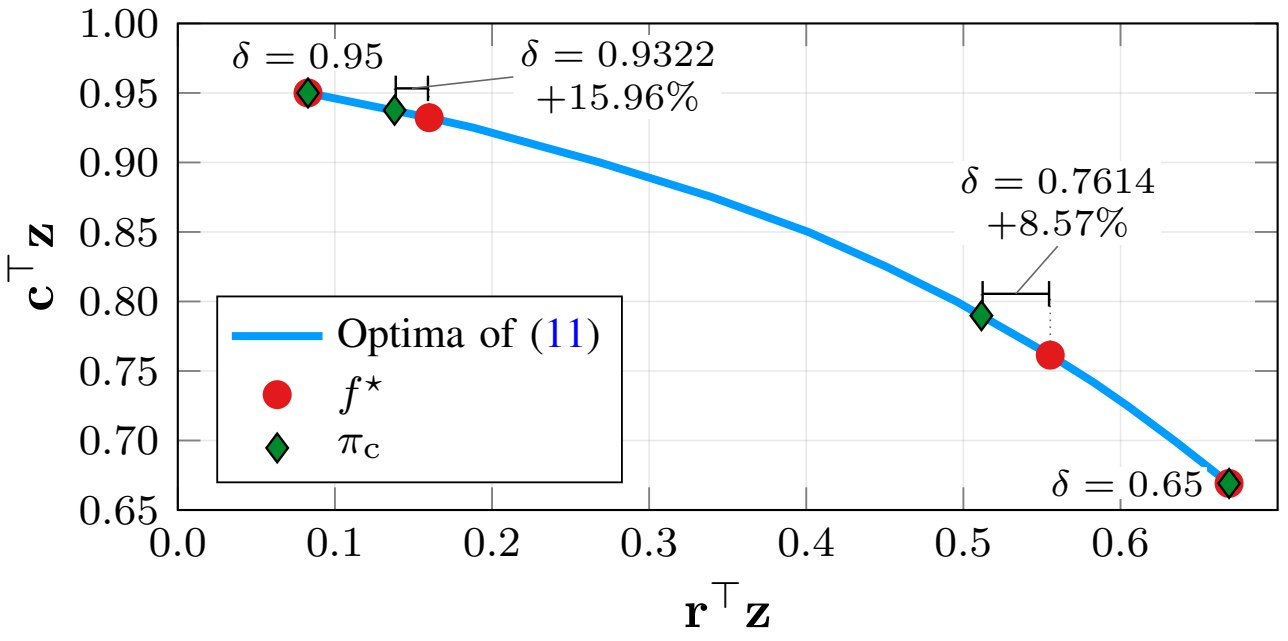

\begin{table}[t]
\centering
\caption{Comparison of reach probability and safety between the optimal switching policy $f\opt$ and the stationary policy $\pi_c$ at selected safety thresholds.}
\label{tab:zermelo_policy_comparison}
\scriptsize
\setlength{\tabcolsep}{3pt}
\begin{tabular}{c|cc|cc|c}
\toprule
\multirow{2}{*}{$\delta$}
& \multicolumn{2}{c|}{LP~\eqref{eq:lp_cons} and $f\opt$}
& \multicolumn{2}{c|}{$\pi_c$}
& \multirow{2}{*}{\shortstack{Reach Probability\\Increase (\%)}} \\
\cmidrule(lr){2-3}
\cmidrule(lr){4-5}
& $\bm r\tr\bm z \!=\! g_{\mathrm r}^{f\opt}\!(\alpha)$
& $\bm c\tr\bm z \!=\! g_{\mathrm c}^{f\opt}\!(\alpha)$
& $g_{\mathrm r}^{\pi_c}\!(\alpha)$
& $g_{\mathrm c}^{\pi_c}\!(\alpha)$
& \\
\midrule
$0.6500$ & $0.6691$ & $0.6691$ & $0.6691$ & $0.6691$ & $0.00$ \\
$0.6692$ & $0.6690$ & $0.6692$ & $0.6204$ & $0.7122$ & $7.83$ \\
$0.6750$ & $0.6624$ & $0.6750$ & $0.6204$ & $0.7122$ & $6.78$ \\
$0.7420$ & $0.5826$ & $0.7420$ & $0.5575$ & $0.7600$ & $4.51$ \\
$0.7614$ & $0.5553$ & $0.7614$ & $0.5115$ & $0.7898$ & $8.57$ \\
$0.7750$ & $0.5343$ & $0.7750$ & $0.5115$ & $0.7898$ & $4.46$ \\
$0.8250$ & $0.4510$ & $0.8250$ & $0.4354$ & $0.8334$ & $3.58$ \\
$0.8750$ & $0.3402$ & $0.8750$ & $0.3402$ & $0.8750$ & $0.00$ \\
$0.9000$ & $0.2686$ & $0.9000$ & $0.2678$ & $0.9002$ & $0.27$ \\
$0.9250$ & $0.1881$ & $0.9250$ & $0.1822$ & $0.9266$ & $3.23$ \\
$0.9322$ & $0.1600$ & $0.9322$ & $0.1380$ & $0.9375$ & $15.96$ \\
$0.9500$ & $0.0829$ & $0.9500$ & $0.0829$ & $0.9500$ & $0.00$ \\
\bottomrule
\end{tabular}
\end{table}

\Cref{fig:lp_ZN} and \cref{tab:zermelo_policy_comparison}
illustrate the safety--performance trade-off for
LP~\eqref{eq:lp_cons} and the policies constructed
using~\eqref{eq:policy-construction}
and~\eqref{eq:opt_policy_krass}.
Both policies are constructed from the same optimal LP
solutions, with results reported at 12 thresholds
within $[0.65,0.95]$.
The policy $f\opt$ is evaluated independently on the
augmented state space $\mathcal S\times\{0,1\}$.
Its gains agree with the LP values to the reported
precision and are therefore reported jointly with
the LP values in~\cref{tab:zermelo_policy_comparison}.

At $\delta=0.65$, the safety constraint is inactive,
and both policies achieve a reach probability and safety
gain of $0.6691$.
At some thresholds, the stationary policy
$\pi_{\mathrm c}$ is feasible but more conservative,
yielding a smaller reach probability.
For example, at $\delta=0.9322$, $f\opt$ achieves
a reach probability of $0.1600$ and safety gain of
$0.9322$, whereas $\pi_{\mathrm c}$ achieves
$0.1380$ and $0.9375$, respectively.
Where the policies differ, the relative increase in reach
probability over $\pi_{\mathrm c}$ ranges from
$0.27\%$ to $15.96\%$, computed from unrounded gains.

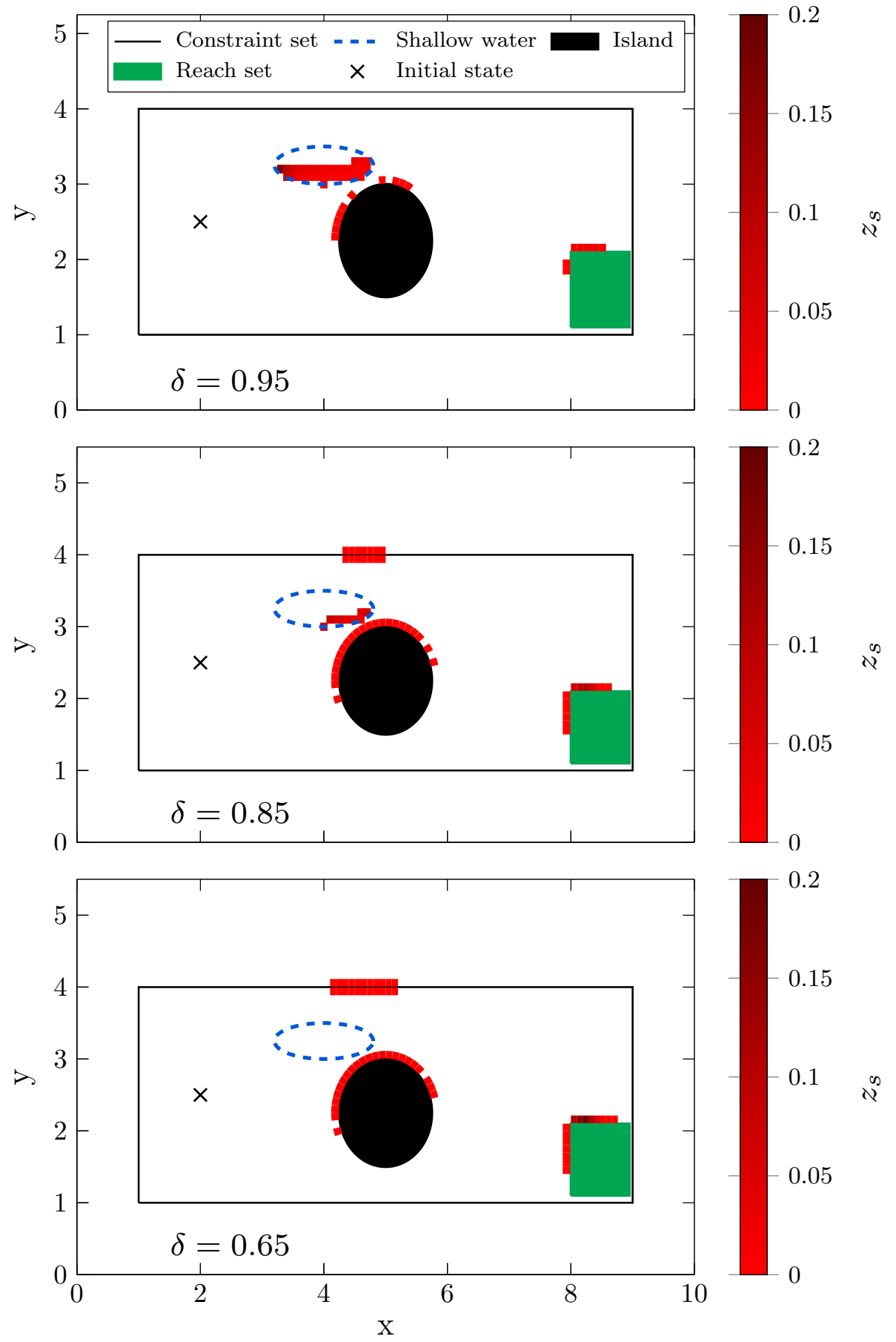
\begin{figure}[t]
    \centering
    \begin{tabular}{@{}c@{}}
        \begin{tikzpicture}
            \node[anchor=south west, inner sep=0] (panel) at (0,0)
                {\includegraphics[height=4.5cm]
                {\detokenize{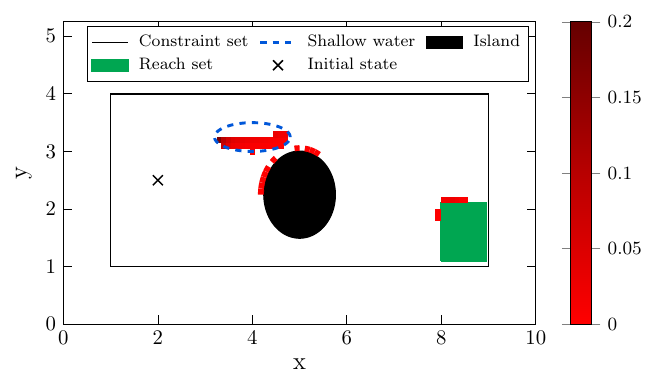}}};

            \begin{scope}[
                x={(panel.south east)},
                y={(panel.north west)}
            ]
                \node[anchor=west, font=\footnotesize,
                      fill=white, inner sep=1pt]
                    at (0.20,0.22) {$\delta=0.95$};

                \node[rotate=90, font=\footnotesize,
                      fill=white, inner sep=1pt]
                    at (1.02,0.53) {$z_s$};
            \end{scope}
        \end{tikzpicture}
        \\[-7.7mm]

        \begin{tikzpicture}
            \node[anchor=south west, inner sep=0] (panel) at (0,0)
                {\includegraphics[height=4.5cm]
                {\detokenize{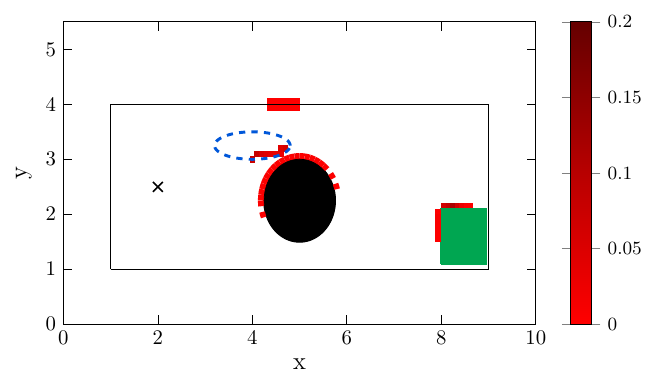}}};

            \begin{scope}[
                x={(panel.south east)},
                y={(panel.north west)}
            ]
                \node[anchor=west, font=\footnotesize,
                      fill=white, inner sep=1pt]
                    at (0.20,0.22) {$\delta=0.85$};

                \node[rotate=90, font=\footnotesize,
                      fill=white, inner sep=1pt]
                    at (1.02,0.53) {$z_s$};
            \end{scope}
        \end{tikzpicture}
        \\[-7.7mm]

        \begin{tikzpicture}
            \node[anchor=south west, inner sep=0] (panel) at (0,0)
                {\includegraphics[height=4.5cm]
                {\detokenize{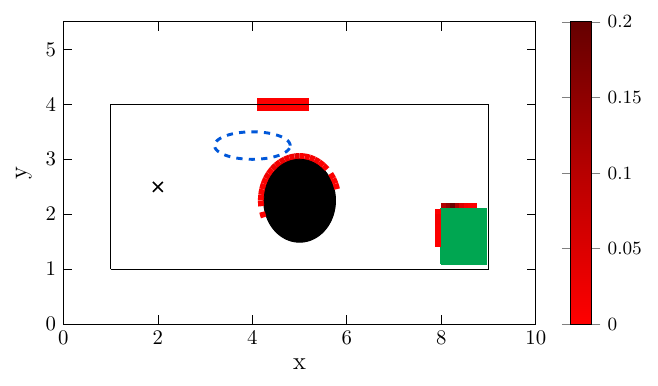}}};

            \begin{scope}[
                x={(panel.south east)},
                y={(panel.north west)}
            ]
                \node[anchor=west, font=\footnotesize,
                      fill=white, inner sep=1pt]
                    at (0.20,0.22) {$\delta=0.65$};

                \node[rotate=90, font=\footnotesize,
                      fill=white, inner sep=1pt]
                    at (1.02,0.53) {$z_s$};
            \end{scope}
        \end{tikzpicture}
        \\[-3mm]
    \end{tabular}

        \caption{Long-run state occupation $z_s$ under $f\opt$
                 for $\delta=0.95$, $0.85$, and $0.65$ (top to bottom);
                 darker red indicates larger mass.}
    \label{fig:krass-policy-distributions}
\end{figure}

\begin{table}[t]
    \centering
\caption{Outcome probabilities under $f\opt$:
remaining safe forever ($\mathbb{P}_{\mathrm{safe}}$),
port arrival ($\mathbb{P}_{\mathrm{port}}$),
permanent shallow-water anchoring ($\mathbb{P}_{\mathrm{sw}}$),
and avoid-set entry ($\mathbb{P}_{\mathrm{avoid}}$).}
    \label{tab:krass-policy-summary}
    \scriptsize
    \setlength{\tabcolsep}{4pt}
    \begin{tabular}{c c c c c}
        \toprule
        $\delta$
        & $\mathbb{P}_{\mathrm{safe}}$
        & $\mathbb{P}_{\mathrm{port}}$
        & $\mathbb{P}_{\mathrm{sw}}$
        & $\mathbb{P}_{\mathrm{avoid}}$ \\
        \midrule
        $0.95$ & $0.9500$ & $0.0829$ & $0.8671$ & $0.0500$ \\
        $0.85$ & $0.8500$ & $0.4013$ & $0.4487$ & $0.1500$ \\
        $0.65$ & $0.6691$ & $0.6691$ & $0.0000$ & $0.3309$ \\
        \bottomrule
    \end{tabular}
\end{table}
Figure \ref{fig:krass-policy-distributions} presents the long-run state occupation of $f\opt$ for safety levels $\delta=0.95$ (top), $\delta=0.85$ (middle), and $\delta=0.65$ (bottom). The red color represents the occupation mass $z_s$ in the shallow-water region and the absorption mass on the reach and avoid boundaries. Darker red corresponds to higher values. Boundary mass reflects the probability of entry at the boundary, whereas shallow-water occupation indicates permanent anchoring. Outcome probabilities for the three safety levels are summarized in Table \ref{tab:krass-policy-summary}.

In the top plot ($\delta=0.95$), the policy uses anchoring
in shallow water to achieve the required safety level.
As shown in~\cref{tab:krass-policy-summary}, the ship
anchors permanently with probability $0.8671$ and reaches
the port with probability $0.0829$. With the remaining
probability $0.05$, the ship hits the island or crosses
the river boundary. Reducing the safety threshold to $\delta=0.85$ (middle plot)
increases the probability of reaching the port to $0.4013$,
while the probability of permanent anchoring in shallow
water decreases to $0.4487$. The ship accepts a higher risk
of hitting the island or crossing the river boundary,
with this probability increasing to $0.15$.
Finally, at $\delta=0.65$ (bottom plot), the safety constraint
is inactive and permanent anchoring has zero probability.
The ship reaches the port with probability $0.6691$ or
hits the island or crosses the river boundary with
probability $0.3309$. Both the constructed stationary
policy $\pi_c$ and the switching policy $f\opt$ attain
the LP optimum to numerical precision.

\section{Conclusion} \label{sec:con}
This study addresses the probabilistic reach--avoid problem
by formulating it as a constrained multichain average-reward
MDP with absorbing reach and avoid sets. The reward gain
represents the probability of reaching the reach set before
entering the avoid set, while the safety gain denotes the
probability of never entering the avoid set. These gain
identities lead to LP~\eqref{eq:lp_cons}, which characterizes
the trade-off between safety and performance. Although
stationary randomized policies may fail to attain the LP
optimum, a lottery-based construction using two stationary
policies and a single randomized switch yields an optimal
policy whenever the LP is feasible. For pure avoidance,
superlevel sets of the optimal safety gain characterize
the probabilistic safe sets. The inverted pendulum and
Zermelo navigation examples demonstrate these results.
In the Zermelo example, the switching policy increases
the reach probability by up to $15.96\%$ relative to the
constructed stationary policy over the tested safety
thresholds, while satisfying the required safety level.
Future research will address discretization error and
computational scalability.

\bibliographystyle{IEEEtran}
\bibliography{Ref}

\end{document}